\documentclass[11pt]{article}
\usepackage{amsmath,latexsym,amssymb}
\def\picsfolder{.}
\usepackage{url}
\usepackage{lscape}
\usepackage{pdflscape}
\usepackage{rotating}
\usepackage{graphics}
\usepackage{enumerate}
\newif\ifpdf
\pdffalse
\ifpdf
\def\fext{pdf}
\else
\def\fext{eps}
\fi
\DeclareMathAlphabet{\mathBB}{U}{bbold}{m}{n}
\providecommand{\1}{\mathBB{1}}
\newcommand{\ignore}[1]{}
\newcommand{\mfalls}{\quad\mbox{if \;}}

\renewcommand{\cases}[1]{\left\{\begin{array}{rl}#1\end{array}\right.}

\newcommand{\mbu}{\quad\mbox{ and }\quad}
\newcommand{\mbs}[1]{\mbox{ \;#1\; }}

\newcommand{\mbsr}[1]{\mbox{#1\; }}

\newcommand{\mf}{\quad\mbox{\;for \;}}
\newcommand{\mfa}{\quad\mbox{\;for all \;}}

\renewcommand{\P}{\mathbf{P}}

\def\given{\hspace{0.8pt}|\hspace{0.8pt}}

\def\ggiven{\hspace{1.2pt}\Big|\hspace{1.2pt}}

\newcommand{\R}{{\mathbb{R}}}
\newcommand{\N}{{\mathbb{N}}}

\newcommand{\C}{{\mathbb{C}}}
\newcommand{\Z}{{\mathbb{Z}}}

\newcommand{\vp}{\varphi}

\def\hatf{\widehat f}

\newcommand{\bec}{\begin{equation}}
\newcommand{\eec}{\end{equation}}
\newcommand{\bac}{\begin{eqnarray}}
\newcommand{\eac}{\end{eqnarray}}
\newcommand{\be}{\begin{displaymath}}
\newcommand{\ee}{\end{displaymath}}
\newcommand{\ba}{\begin{eqnarray*}}
\newcommand{\ea}{\end{eqnarray*}}

\newcommand{\equ}[1]{(\ref{#1})}

\providecommand{\eop}{{}\hfill {}\hfill{$\Box
$}\vspace{0.3cm}\pagebreak[2]\par}
\def\SPT{\mathsf{ST}}
\newtheorem{proposition}{Proposition}[section]
\newtheorem{theorem}[proposition]{Theorem}

\newtheorem{lemma}[proposition]{Lemma}

\newtheorem{uremark}[proposition]{Remark}
\newenvironment{remark}{\begin{uremark}\em}{\hfill$\Diamond$\end{uremark}}
\newtheorem{uexample}[proposition]{Example}

\newcommand{\Section}[1]{\section{#1}\setcounter{figure}{0}\setcounter{table}{0}\setcounter{equation}{0}}

\def\qed{\mbox{$\Box$}}

\newcommand{\weight}{\mathop{\mbox{\sf weight}}}

\newtheorem{udefi}{Definition}[section]
\newtheorem{utheo}{Theorem}
\newtheorem{usatz}[udefi]{Satz}
\newtheorem{uprop}[udefi]{Proposition}

\newtheorem{ubemerkung}[udefi]{Bemerkung}
\newtheorem{ukorollar}[udefi]{Korollar}

\def\Re{\mathrm{Re}}
\def\Im{\mathrm{Im}}

\newcommand{\rest}[1]{\raisebox{-0.3em}{$\big|_{\scriptstyle
#1}$}}

\begin{document}
\title{The Random Spanning Tree on Ladder-like Graphs}
\author{
Achim Klenke\\
Institut f\"{u}r Mathematik\\
Johannes Gutenberg-Universit\"{a}t Mainz\\
Staudingerweg 9\\
55099 Mainz}
\date{\small  28.03.2017}
\maketitle
\begin{abstract}
Random spanning trees are among the most prominent determinantal point processes. We give four examples of random spanning trees on ladder-like graphs whose rungs form stationary renewal processes or regenerative processes of order two, respectively. Up to a trivial thinning with additional coin flips, for each of the first two examples the renewal processes exhaust the whole class of stationary regenerative (of order one) determinantal point processes. We also give an example of a regenerative process of order two that has no representation in terms of a random spanning tree.

Our examples illustrate a theorem of Lyons and Steif (2003) which characterizes regenerative determinantal point processes in terms of their Fourier transform.  For the regenerative process, we also establish a Markov chain description in the spirit of H\"{a}ggstr\"{o}m (1994).

On the technical side, a systematic counting scheme for random spanning trees is developed that allows to compute explicitly the probabilities. In some cases an electrical network point of view simplifies matters. \end{abstract}
\Section{Introduction}
\label{S1}
\subsection{Random Spanning Trees}
\label{S1.1}
On a finite connected (undirected) graph $G=(V,E)$ with vertex set $V$ and edge set $E$, there is a finite set $\SPT(G)\subset 2^E$ of spanning trees $t$; that is, the graph $(V,t)$ is connected and loop-free (see, e.g., \cite{Bollobas1998}). The uniform distribution on $\SPT(G)$ is called the uniform spanning tree measure. If we assign to each edge $e\in E$ a weight $\weight(e)\in[0,\infty)$ and define $\weight(t)=\prod_{e\in t}\weight(e)$ and $\weight(F)=\sum_{g\in F}\weight(g)$ for $F\subset 2^E$, then $\P[F]:=\weight(F)/\weight(\SPT(G))$, $F\subset\SPT(G)$, defines the weighted spanning tree measure with weight function $\weight$.

For an infinite graph $G=(V,E)$, in order to define a uniform or weighted spanning tree measure, one can exhaust $G$ by an increasing sequence of finite subgraphs $G_n$ and define the measure $\P$ as the limit of the uniform (or weighted) spanning tree measures $\P_n$ on $G_n$. By a simple monotonicity argument, it is shown that this limit exists (see, e.g., \cite[Proposition 5.6]{BenjaminiLyonsPeresSchramm2001}) and does not depend on the choice of the sequence $(G_n)$. $\P$
is called the \emph{free spanning forest measure}. The reason for this name is that, in general, the limit $\P$ is not necessarily concentrated on connected graphs (but clearly on loop-free spanning graphs); that is, on spanning forests. For $G=\Z^d$ the $d$-dimensional integer lattice, it is shown in \cite{Pemantle1991} that the uniform spanning forest measure is indeed concentrated on trees if and only if $d\leq4$. On the other hand, for $d\geq 5$, the free spanning forest measure is concentrated on forests with infinitely many trees.

A simple argument that uses Wilson's method of simulating random spanning trees shows that $\P$ is concentrated on trees if the graph random walk with transition probabilities proportional to the edge weights is recurrent (\cite[Proposition 5.6]{BenjaminiLyonsPeresSchramm2001}). The graphs that we consider in this paper are essentially one-dimensional. Hence recurrence will be obvious and $\P$ can be considered as the uniform (or weighted) spanning tree measure. In the sequel, $T$ will always be the generic random spanning tree under the measure $\P$.

There is a rich body of literature on random spanning trees, see, e.g.
\cite{AizenmanBurchardNewman1997,AldousSteele1992,Aldous1990,BenjaminiLyonsPeresSchramm2001,
BurtonPemantle1993,
Haggstrom1994,Haggstrom1995b,LawlerSchrammWerner2004,LyonsPeres2016,
Pemantle1991,ProppWilson1998b,Wilson1996}. Particularly interesting is a connection to electrical networks that was first discovered by Gustav Kirchhoff in 1847. In order to describe Kirchhoff's results, consider our (finite) graph as an electrical network with nodes $v\in V$ and resistors with resistances $R(e)=1/\weight(e)$ along the edges $e\in E$. Denote by $R_{\rm eff}(v,w)$ the effective resistance between the nodes $v,w\in V$. Kirchhoff proved in \cite{Kirchhoff1847} that for each edge $e=\{v,w\}$, we have
\begin{equation}
\label{E1.01}
\P[e\in T]=\frac{R_{\rm eff}(v,w)}{R(e)}.
\end{equation}

To put the result a bit differently, assume that we hook a battery at $v$ and $w$ such that a unit current flows into the network at $v$ and out of the network at $w$. Denote by $I_{e,f}$ the current along the edge $f$ (where we assume that we have given all edges in $E$ an arbitrary but fixed orientation in order to fix the sign of the current). Then by Ohm's law, we have $I_{e,e}=\frac{R_{\rm eff}(v,w)}{R(e)}$, that is $\P[e\in T]=I_{e,e}$.

Burton and Pemantle's transfer current theorem \cite{BurtonPemantle1993} describes the full random spanning tree (also in infinite networks) in terms of the electrical network: Let $n\in\N$ and let $e_1,\ldots,e_n\in E$ be pairwise distinct edges. Then the probability that all these edges are in the random spanning tree can be expressed in terms of the determinant of a submatrix of the so-called impedance (or transfer current) matrix $(I_{e,f})_{e,f\in E}$,
\begin{equation}
\label{E1.02}
\P[e_1,\ldots,e_n\in T]=\det\Big((I_{e_i,e_j})_{i,j=1,\ldots,n}\Big).
\end{equation}
Note that the signs of the individual entries $I_{e,f}$ depend on the arbitrary choice of the orientations of the edges but the value of the determinant does not.
More information about electrical networks can be found, e.g., in \cite{DoyleSnell1984} and \cite{LyonsPeres2016}. A streamlined proof of the transfer current theorem \equ{E1.02} is given in \cite{BenjaminiLyonsPeresSchramm2001}.
A nice introduction to random spanning trees and electrical networks can be found in Jarai's lecture notes \cite{Jarai2009}.

\subsection{Determinantal Point Processes}
\label{S1.2}
Equation \equ{E1.02} states that $T$ is a \emph{determinantal point process} on the index set $E$. More generally, if $R$ is a countable set and $A$ is a Hermitian matrix indexed by $R$, a $\{0,1\}^R$-valued process $X$ is called a determinantal point process if for any finite subset $S\subset R$, we have
\begin{equation}
\label{E1.03}
\P\big[X(s)=1\mbs{for all}s\in S\big]=\det\big((A_{ij})_{i,j\in S}\big).
\end{equation}

A sufficient condition for $A$ to be the matrix of some determinantal point process is that it is a positive contraction (see~\cite[Theorem 8.1]{Lyons2003} or \cite[Theorem~1.1]{ShiraiTakahashi2003II} for real-valued $A$). By a theorem of Soshnikov (see \cite[Theorem 3]{Soshnikov2000}), if $A$ is a Hermitian operator on $\ell^2(R)$, then $A$ is the matrix of a determinantal point process if and only if $A$ and $I-A$ are nonnegative definite (where $I$ is the identity matrix). The necessity of this condition is clear: Equation \equ{E1.03} implies that $A$ is nonnegative definite. On the other hand, a simple application of the inclusion/exclusion formula shows that the process of zeros, that is $X'(s)=1-X(s)$, fulfills \equ{E1.03} with $A'=I-A$. Hence $X'=1-X$ is also a determinantal point process and has the matrix $I-A$. For details, see, e.g. \cite[Section 2]{Gottlieb2007}.

For the general theory, the matrix $A$ need not even be Hermitian, but in this article we restrict ourselves to the Hermitian case. In order to clarify this in the theorems, we will refer to $X$ as a Hermitian determinantal point process.

In \cite{LyonsSteif2003} Lyons and Steif investigate \emph{stationary} determinantal point processes indexed by the integer lattice $\Z$. That is, $A$ is Hermitian and fulfills $A_{k,l}=A_{0,l-k}$ for all $k,l\in\Z$. In other words, $A$ is an infinite Hermitian Toeplitz matrix. Taking Fourier transforms, there exists a (unique for $A$) measurable function $f:[0,1]\to[0,1]$ such that $A_{0,k}=\hatf(k)$ for all $k\in\Z$, where (with $\lambda$ the Lebesgue measure)
\begin{equation}
\label{E1.04}
\hatf(k)=\int_{[0,1]}f(x)\,e^{-2\pi i\,kx}\,\lambda(dx).
\end{equation}
In fact, by Bochner's theorem, for the nonnegative definite matrix $A$, there exists a measure $\mu$ on $[0,1)$ such that $A_{0,k}=\int e^{-2\pi i\,kx}\,\mu(dx)$. Since also $I-A$ is nonnegative definite, there also exists a measure $\mu'$ on $[0,1)$ such that $I_{0,k}-A_{0,k}=\int e^{-2\pi i\,kx}\,\mu'(dx).$
Summing up, we get
$I_{0,k}=\int e^{-2\pi i\,kx}\,(\mu+\mu')(dx)$
and hence $\mu+\mu'=\lambda$. Now define $f=d\mu/d\lambda$.
By Fourier inversion, we have
\begin{equation}
\label{E1.05}
f(x)=\sum_{k\in\Z}A_{0,k}\,e^{i\,kx}.
\end{equation}
The stationary process $X$ is said to be a \emph{regenerative process of order $k\in\N$} if it renews after $k$ successive $1$s. More precisely, for measurable $B\subset\{0,1\}^{\{-k,-k-1,\ldots\}}$ and $C\subset\{0,1\}^\N$, we have
$$\begin{aligned}
\P\Big[X\rest{\N}\in C\ggiven X\rest{\{-k,-k-1,\ldots\}}&\in B,\,X(0)=X(-1)=\ldots=X(-k+1)=1\Big]\\
&=\P\Big[X\rest\N\in C\ggiven X(0)=X(-1)=\ldots=X(-k+1)=1\Big].\end{aligned}
$$
Lyons and Steif \cite[Proposition 2.10]{LyonsSteif2003} give a simple characterization of renerative determinantal point processes in terms of the Fourier transform $f$ of $A$:
\begin{proposition}[Lyons and Steif, Proposition 2.10 of \cite{LyonsSteif2003}]
\label{P1.1}
Let $X$ be a stationary determinantal point process with Hermitian matrix $A$ and Fourier transform $f$. Let $k\in\N$. The process $X$ is regenerative of order $k$ if and only if
$1/f$ is a trigonometric polynomial of order at most $k$, that is, there exist numbers $c_0,\ldots,c_k,\varphi_1,\ldots,\varphi_k\in\R$ such that
\begin{equation}
\label{E1.06}
\frac{1}{f(x)}=c_0+\sum_{j=1}^kc_j\,\cos\big(2\pi j(x+\varphi_j)\big),\qquad x\in[0,1].
\end{equation}
\end{proposition}
Lyons and Steif give an example (due to Soshnikov) for $k=1$, the case of a classical renewal process: \begin{equation}
\label{E1.07}
\hatf(k)=\frac{1-\alpha}{1+\alpha}\,\alpha^{|k|},\mfa k\in\Z,
\end{equation}
where $\alpha\in(0,1)$ is a parameter of the model. In this case
$$f(x)=\frac{(1-\alpha)^2}{1+\alpha^2-2\alpha\,\cos(2\pi x)}.$$
In fact, the renewal distribution (waiting time for the next $1$) is given by
\begin{equation}
\label{E1.08}
\P\big[\inf\{k\geq1:X_k=1\}=m\given X_0=1\big]=r_m:=
(1-\alpha)^2\,m\,\alpha^{m-1}.
\end{equation}
This waiting time is distributed as the sum of two geometric random variables and we will see later an intuitive reason for this fact.

\begin{remark}
\label{R1.2}
Note that for a given stationary determinantal point process $X$, the Hermitian Toeplitz matrix $A$ is not unique, but for any $\varphi\in\R$, the matrix $A^\varphi_{k,l}:=A_{k,l}\cdot e^{2\pi i\varphi(l-k)}$ describes the same process $X$. If $f$ is the corresponding Fourier transform, then $f^\varphi(x)=f((x+\varphi)\mod 1)$. In particular, if $A$ is an impedance matrix and is hence real-valued, the matrix $A^{1/2}_{k,l}=(-1)^{l-k}A_{k,l}$ describes the same process $X$ and is the impedance matrix that belongs to alternating the orientation of the edges.

If $A$ has only nonzero entries, then it is simple to show that any matrix $A'$ that describes the same process $X$ is of the form $A'=A^\varphi$ for some $\varphi\in\R$. We will discuss this in some more detail in Section~\ref{S4} for a particular matrix $A$.
\end{remark}

\subsection{Main results}
\label{S1.3}
While Soshnikov's example for a stationary determinantal point process remains rather abstract, random spanning trees on ladder-like graphs provide natural and intuitive examples. In this section, we present two graphs for which the weighted random spanning tree gives us Soshnikov's example for all values of $\alpha\in(0,1)$. Furthermore, we present two graphs that are examples for an order two regenerative determinantal point process.

\subsubsection{The simple ladder graph}
\label{S1.3.1}
We start with the simple (two-sided infinite) ladder graph.
Consider the vertex set $V^L=\{0,1\}\times\Z$ and denote by
$$E^L=\big\{z_m,\,h_{0,m},\,h_{1,m}:\;m\in\Z\big\}$$ the set of edges where $$z_m:=\{(0,m),(1,m)\}\mbu h_{i,m}=\{(i,m-1),(i,m)\},\quad m\in\Z,\,i=0,1.$$ The \emph{simple ladder graph} is the graph $G^L:=(V^L,E^L)$. For $m\leq n$, write
$$V^L_{m,n}:=\{0,1\}\times\{m,\ldots,n\}$$
 and denote the induced edge set by $E^L_{m,n}$. Finally, define the finite ladder graph
\begin{equation}
\label{E1.09}
G^L_{m,n}=(V^L_{m,n},E^L_{m,n})
\end{equation}
as the induced subgraph.
\begin{figure}[h]
\label{F1.1}
\begin{picture}(400,110)(0,0)
\put(40,02){\includegraphics[scale=0.8]{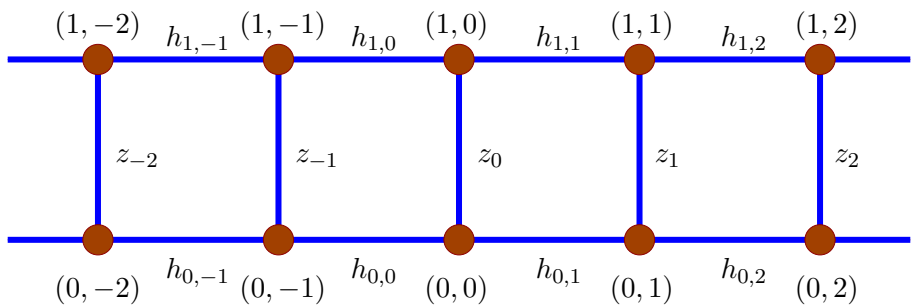}}
\put(200, 00){$(0,0)$}
\put(200,100){$(1,0)$}
\put(270, 00){$(0,1)$}
\put(340,100){$(1,2)$}
\put(340, 00){$(0,2)$}
\put(270,100){$(1,1)$}
\put(130, 00){$(0,-1)$}
\put(130,100){$(1,-1)$}
\put(060, 00){$(0,-2)$}
\put(060,100){$(1,-2)$}
\put(172,7){$h_{0,0}$}
\put(172,95){$h_{1,0}$}
\put(242,7){$h_{0,1}$}
\put(242,95){$h_{1,1}$}
\put(312,7){$h_{0,2}$}
\put(312,95){$h_{1,2}$}
\put(102,7){$h_{0,-1}$}
\put(102,95){$h_{1,-1}$}
\put( 83,51){$z_{-2}$}
\put(151,51){$z_{-1}$}
\put(220,51){$z_{0}$}
\put(287,51){$z_{1}$}
\put(355,51){$z_{2}$}
\end{picture}
\caption[b]{A finite section of the simple ladder graph}
\end{figure}
We will be interested in the weight function $\weight(z_m)=c$ for all $m\in\Z$ for some $c>0$ and $\weight(h_{i,m})=1$ for all $m\in\Z$ and $i=0,1$. Denote by $\P^c_n$ the weighted spanning tree distribution on $\SPT\big(G^L_{-n,n}\big)$, that is
\begin{equation}
\label{E1.10}
\P^c_n[\{t\}]=\frac{\weight(t)}{\weight\big(\SPT\big(G^L_{-n,n}\big)\big)}\mf t\in\SPT\big(G^L_{-n,n}\big).
\end{equation}

The random walk on $V$ that jumps from $(i,n)$ to $(1-i,n)$ with probability $c/(c+2)$ and to $(i,n+1)$ and $(i,n-1)$ each with probability $1/(c+2)$ is clearly recurrent. Hence the limit $\P^c=\lim_{n\to\infty}\P^c_n$ is indeed a probability measure that concentrates on spanning trees of $G^L$. Denote by $T$ the generic (weighted) random spanning tree under $\P^c$ and define $X=(X_m)_{m\in\Z}$ by $X_m=\1_T(z_m)$. For $c=\infty$, we define $X\equiv1$. This makes sense since giving the rungs infinite weight is the same as conditioning on all rungs to be in the spanning tree. Let $p\in[0,1]$ and let $(Y_m)_{m\in\Z}$ be an iid sequence (independent of $X$) of Bernoulli random variables mit parameter $p$. Finally, define $\tilde X_m:=X_mY_m$.

Obviously, $X$ is a stationary renewal process with some real-valued symmetric Toeplitz matrix $A$ and Fourier transform $f$. Flipping extra coins for each $1$ of $X$ does not change the renewal property. Hence $\tilde X$ is also a renewal process. Furthermore, it is determinantal with matrix $\tilde A=pA$.

\begin{theorem}
\label{T1}

Let $\alpha\in[0,1)$ and $p\in[0,1]$. Define
\begin{equation}
\label{E1.11}
c:=\frac{(1-\alpha)^2}{2\alpha}>0 \mfalls \alpha\in(0,1)
\end{equation}
and $c=\infty\;$ if $\;\alpha=0$.
\begin{enumerate}[(i)]
\item
For the simple ladder graph, under the weighted spanning tree measure $\P^c$, the process $X$ is determinantal with matrix
$A_{k,l}=\hatf(l-k)$, where
\begin{equation}
\label{E1.12}
\hatf(k)=\frac{1-\alpha}{1+\alpha}\,\alpha^{|k|} \mfa k\in\Z
\end{equation}
and
\begin{equation}
\label{E1.13}
f(x)=\frac{c}{\displaystyle c+1-\cos(2\pi x)}.
\end{equation}

In particular, $X$ is a renewal process with renewal distribution given by \equ{E1.08}.

\item
Furthermore, the thinned process $\tilde X$ is determinantal with matrix $\tilde A=pA$ and Fourier transform
\begin{equation}
\label{E1.14}
\tilde f(x)=\frac{pc}{\displaystyle c+1-\cos(2\pi x)}.
\end{equation}
\item
Every stationary (Hermitian) determinantal point process that is also a renewal process has the distribution of $\tilde X$ for some $\alpha\in[0,1)$ and $p\in[0,1]$.
\end{enumerate}
\end{theorem}

Note that (ii) is a direct consequence of (i). For (iii), note that by Proposition~\ref{P1.1}, for every such such process, we have $1/\tilde f(x)=c_0+c_1\cos(2\pi (x+\vp))$ for some $c_0,c_1,\vp\in\R$. Without loss of generality we may assume $c_1\leq0$, otherwise take $\varphi+\pi$ instead of $\varphi$. By Remark~\ref{R1.2}, the function $f^{-\varphi}(0)=1/(c_0+c_1\cos(2\pi x))$ describes the same process. Hence, we may assume $\varphi=0$. Recall that $f$ takes values in $[0,1]$, hence $c_0+c_1\geq1$. Letting
$$p:=\frac{1}{c_0+c_1}\mbu c:=-\frac{1}{c_1p},$$
we get \equ{E1.14}.

Note that for $\alpha=0$, we have $A=I$ the identity matrix, that is, $X\equiv1$, resulting in $\tilde f\equiv p$.

\begin{remark}
\label{R1}
The very form of $\hatf$ shows that for $p=1$, the renewal distribution of $X$ is $\delta_1\ast\gamma_\alpha\ast\gamma_\alpha$ where $\gamma_\alpha$ is the geometric distribution on $\N_0$ with parameter $\alpha$. Between two renewal events at times $m<n$ say, there is exactly one horizontal edge $h_{i,k}$, $m<k\leq n$, $i=0,1$ missing. Since the very position of the missing edge does not change the weight of the tree, each choice has the same probability. This allows for a very simple construction of the spanning tree. Assuming $z_0\in T$, let $Y,Y'$ be independent and $\gamma_\alpha$ distributed random variables. Furthermore let $W$ be an independent Bernoulli 1/2 random variable. Given $Y+Y'$, the random variable $Y$ is uniformly distributed on $\{0,\ldots,Y+Y'\}$. Hence $Y+Y'+1$ marks the next renewal event and the missing horizontal edge is $h_{W,Y+1}$. For the next renewal interval proceed similarly.

In order to get rid of the assumption $z_0\in T$, the first renewal time has to be chosen in a size-biased manner.
\end{remark}

In Section~\ref{S2}, we will give different approaches with explicit calculations for the matrix $A$. The heart of the computations is an explicit formula for the weighted number of spanning trees on a finite part $G^L_{0,n-1}$ (recall \equ{E1.09}) of the ladder graph.

\begin{proposition}
\label{P1.3}
Assume that $\weight(h_{i,n})=1$ and $\weight(z_n)=c\in(0,\infty)$ for all $i=0,1$, $n\in\Z$. Let
\begin{equation}
\label{E1.15}
\alpha:=c+1-\sqrt{c^2+2c}\;\in\,(0,1).
\end{equation}
The weighted number of spanning trees on $G^L_{0,n-1}$ is
\begin{equation}
\label{E1.16}
\weight\big(\SPT\big(G^L_{0,n-1}\big)\big)=\frac12\frac{1-\alpha}{1+\alpha}\big(\alpha^{-n}-\alpha^{n}\big).
\end{equation}

In particular, for $c=1$, the number of spanning trees on $G^L_{0,n-1}$ is
$$\#\SPT\big(G^L_{0,n-1}\big)=\frac{1}{2\sqrt{3}}\left(\big(2+\sqrt{3}\,\big)^n-\big(2-\sqrt{3}\,\big)^n\right).$$
\end{proposition}

\subsubsection{The zigzag ladder graph $G^{H_2}$}
\label{S1.3.2}
A different example for a random spanning tree that yields a renewal process is the zigzag graph. It is quite similar to the simple ladder graph but also has some similarities with the helix-3-graph that will be presented later.

\begin{figure}[ht]
\begin{picture}(400,120)(00,0)
\put(20,02){\includegraphics[scale=0.9]{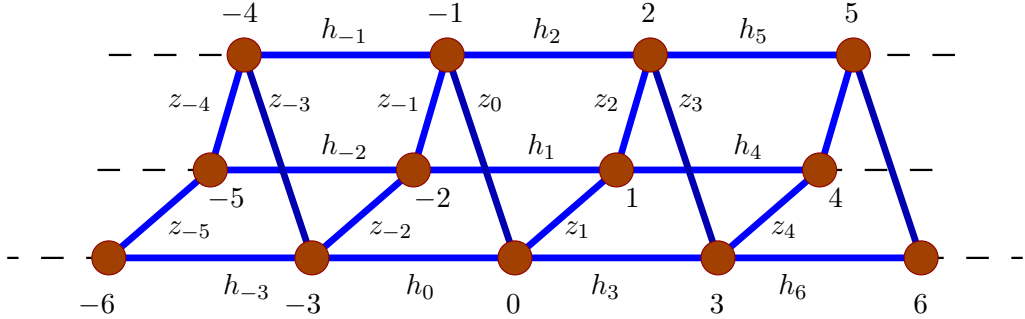}}
\put(048,00){$-6$}
\put(126,00){$-3$}
\put(210,00){$0$}
\put(287,00){$3$}
\put(364,00){$6$}
\put(97,40){$-5$}
\put(255,40){$1$}
\put(332,40){$4$}
\put(175,40){$-2$}
\put(102,110){$-4$}
\put(180,110){$-1$}
\put(261,110){$2$}
\put(338,110){$5$}
\put(313,7){$h_6$}
\put(103,7){$h_{-3}$}
\put(172,7){$h_{0}$}
\put(242,7){$h_{3}$}
\put(140,60){$h_{-2}$}
\put(218,60){$h_{1}$}
\put(296,60){$h_{4}$}
\put(140,104){$h_{-1}$}
\put(220,104){$h_{2}$}
\put(298,104){$h_{5}$}
\put( 82,30){$z_{-5}$}
\put(158,30){$z_{-2}$}
\put(232,30){$z_{1}$}
\put(310,30){$z_{4}$}
\put(120,78){$z_{-3}$}
\put(199,78){$z_{0}$}
\put(275,78){$z_{3}$}
\put(82,78){$z_{-4}$}
\put(161,78){$z_{-1}$}
\put(243,78){$z_{2}$}
\end{picture}
\caption[b]{A finite section of the helix-3-graph}
\label{F1.3}
\end{figure}

\begin{figure}[ht]
\label{F4.2}
\begin{picture}(400,120)(00,0)
\put(20,02){\includegraphics[scale=0.9]{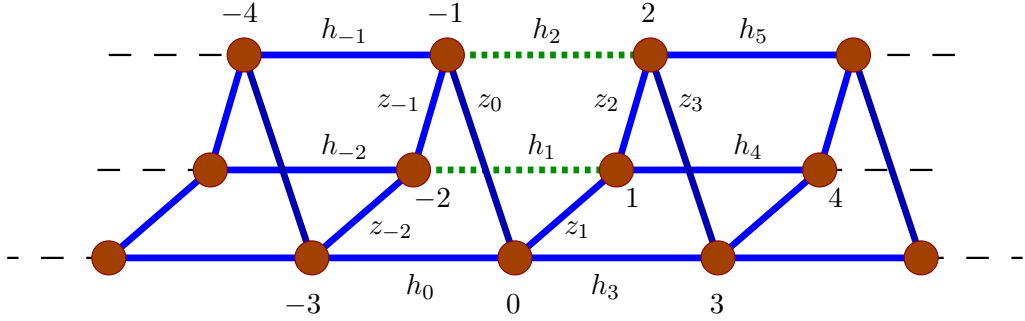}}
\put(126,00){$-3$}
\put(210,00){$0$}
\put(287,00){$3$}
\put(255,40){$1$}
\put(332,40){$4$}
\put(175,40){$-2$}
\put(102,110){$-4$}
\put(180,110){$-1$}
\put(261,110){$2$}
\put(172,7){$h_{0}$}
\put(242,7){$h_{3}$}
\put(140,60){$h_{-2}$}
\put(218,60){$h_{1}$}
\put(296,60){$h_{4}$}
\put(140,104){$h_{-1}$}
\put(220,104){$h_{2}$}
\put(298,104){$h_{5}$}
\put(158,30){$z_{-2}$}
\put(232,30){$z_{1}$}
\put(199,78){$z_{0}$}
\put(275,78){$z_{3}$}
\put(161,78){$z_{-1}$}
\put(243,78){$z_{2}$}
\end{picture}
\caption[b]{Connecting left and right part in the helix-3-graph.}
\end{figure}

Let $k=2,3,\ldots$ and define $V^{H_k}=\Z$ and $E^{H_k}=\{h_m,\,z_m:\,m\in\Z\}$ where $z_m=\{m-1,m\}$ and $h_m=\{m-k,m\}$. Define the helix-$k$-graph $G^{H_k}=(V^{H_k},E^{H_k})$. See Figure~\ref{F4.2} for the helix-2-graph and Figure~\ref{F1.3} below for the helix-3-graph. For $m\leq n$, let $V^{H_k}_{m,n}=\{m,\ldots,n\}$ and define the induced edge set $E^{H_k}_{m,n}$ and the induced subgraph $G^{H_k}_{m,n}$.

Here, we consider the helix-2-graph $G^{H_2}$ that we also call the \emph{zigzag ladder} graph.

\begin{figure}[ht]
\label{F1.2}
\begin{picture}(400,110)(0,0)
\put(63,02){\includegraphics[scale=0.8]{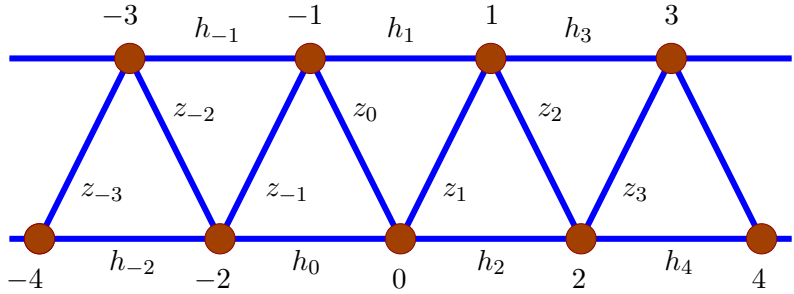}}
\put(064,00){$-4$}
\put(135,00){$-2$}
\put(210,00){$0$}
\put(278,00){$2$}
\put(346,00){$4$}
\put(100,100){$-3$}
\put(170,100){$-1$}
\put(245,100){$1$}
\put(313,100){$3$}
\put(313,7){$h_4$}
\put(103,7){$h_{-2}$}
\put(172,7){$h_{0}$}
\put(242,7){$h_{2}$}
\put(135,95){$h_{-1}$}
\put(208,95){$h_{1}$}
\put(275,95){$h_{3}$}
\put( 92,35){$z_{-3}$}
\put(162,35){$z_{-1}$}
\put(229,35){$z_{1}$}
\put(297,35){$z_{3}$}
\put(127,65){$z_{-2}$}
\put(195,65){$z_{0}$}
\put(265,65){$z_{2}$}
\end{picture}
\caption[b]{A finite section of the zigzag ladder graph (helix-2-graph).}
\end{figure}
Let $c>0$ and consider the weight function $\weight(h_m)=1$ and $\weight(z_m)=c$ for all $m\in\Z$. Note that the corresponding random walk is recurrent. Hence we can define the weighted random spanning tree measure $\P^c$ on $G^{H_2}$ as the limit of $\P^c_n$ on $G^{H_2}_{-n,n}$ as for the simple ladder graph. Let $T$ denote the generic random spanning tree and define $X_m=\1_T(z_m)$.
\begin{theorem}
\label{T2}
With
\begin{equation}
\label{E1.17}
\alpha:=1+\frac{c}{2}-\frac{\sqrt{4c+c^2}}{2}
\end{equation}
or, equivalently, $c=(1-\alpha)^2/\alpha$,
$$f(x)=\frac{c}{\displaystyle c+2-2\cos(2\pi x)}$$
and $\tilde f$ from \equ{E1.12}, the statements of Theorem~\ref{T1} also hold for the zigzag ladder graph.
\end{theorem}
Again the weighted number of spanning trees can be computed explicitly.
\begin{proposition}\label{P1.4}
For $c\in(0,\infty)$, the weighted number of spanning trees on $G^{H_2}_{0,n-1}$ is
$$
\weight\big(\SPT\big(G^{H_2}_{0,n-1}\big)\big)=\frac{1-\alpha}{1+\alpha}\big(\alpha^{1-n}-\alpha^{n-1}\big).
$$
In particular, for $c=1$, the number of spanning trees on $G^{H_2}_{0,n-1}$ is
\begin{equation}
\label{E1.18}
\#\SPT\big(G^{H_2}_{0,n-1}\big)
=\frac{1}{\sqrt{5}}
\left(\left(\frac{3+\sqrt{5}}{2}\right)^{n-1}
-\left(\frac{3-\sqrt{5}}{2}\right)^{n-1}\right).
\end{equation}
\end{proposition}
\subsubsection{The helix-3-graph $G^{H_3}$}
\label{S1.3.3}
In the previous two examples, the process $X$ that resulted from the rungs in the random spanning tree was a renewal process. Here we come to an example where $X$ is regenerative of order $2$ and where we can compute the matrix and the Fourier transform explicitly.

We consider only unit weights on the edges and define the uniform spanning tree measure $\P$ as the limit of the uniform distributions $\P_n$ on $\SPT\big(G^{H_k}_{-n,n}\big)$ (for all $k\geq2$). Let $T$ be the generic uniform spanning tree and let $X_m:=\1_T(z_m)$, $m\in\Z$. Clearly, $X$ is a (Hermitian) determinantal point process that is regenerative of order $k-1$. Hence, by Proposition~\ref{P1.1}, we get that the reciprocal Fourier transform of the matrix $A$ of $X$ is a trigonometric polynomial of degree $k-1$. The main goal of this section is to compute this polynomial explicitly for $k=3$. For larger $k$, it seems to be hopeless to compute the explicit Fourier transform. We will encounter various powers of the golden ratio and hence, for convenience, give it a symbol
\begin{equation}
\label{E1.19}
\gamma:=\frac{\sqrt{5}+1}{2}.
\end{equation}

\begin{theorem}
\label{T3}
\begin{enumerate}[(i)]
\item
For the uniform random spanning tree on the helix-3-graph, the process $X$ is a determinantal point process with (real-valued symmetric Toeplitz) matrix $A_{k,l}=\hatf(l-k)$ given by
\begin{equation}
\label{E1.20}
\hatf(m)=\eta\, \alpha^{|m|}+\bar\eta\,\bar\alpha^{|m|},\qquad m\in\Z,
\end{equation}
where (with $i=\sqrt{-1}$)
\begin{equation}
\label{E1.21}
\eta=\frac{\gamma^{3/2}}{4\sqrt{5}}\;+\;i\,\frac{\gamma^{-3/2}}{4\sqrt{5}}\mbu
\alpha=\frac{\gamma^{-3/2}-1}{2}\;+\;i\,\frac{\gamma^{3/2}-1}{2}.
\end{equation}
\item
The Fourier transform $f$ of $\hatf$ equals
\begin{equation}
\label{E1.22}
f(x)=\frac{1}{4+4\cos(2\pi x)+2\cos(4\pi x)}.
\end{equation}
\end{enumerate}
\end{theorem}
Counting the number of spanning trees with certain additional properties will be an essential tool in the analysis of this example. In order to give a flavour, we present here the explicit formula for the number of spanning trees on a finite subgraph of the helix-3-graph.
\begin{proposition}
\label{P1.5}
For $n\geq3$, the number of spanning trees on $G^{H_3}_{0,n-1}$ is
$$\begin{aligned}
\#\SPT\big(G^{H_3}_{0,n-1}\big)&=\frac{\gamma^2-\gamma^{3/2}}{4\sqrt{5}}\big(\gamma+\sqrt{\gamma}\big)^n
-\frac12
+\frac{\gamma^2+\gamma^{3/2}}{4\sqrt{5}}\big(\gamma-\sqrt{\gamma}\big)^n\\
&\quad
+\frac{-\gamma^{-2}-i\gamma^{-3/2}}{4\sqrt{5}}\big(-\gamma^{-1}+i/\sqrt{\gamma}\big)^n
+\frac{-\gamma^{-2}+i\gamma^{-3/2}}{4\sqrt{5}}\big(-\gamma^{-1}-i/\sqrt{\gamma}\big)^n.
\end{aligned}
$$
\end{proposition}
Note that only the third summand vanishes asymptotically since $\gamma-\sqrt\gamma\approx0.346$. The fourth and fifth summand are complex rotations. More precisely,
$$\limsup_{n\to\infty}\Big|\#\SPT\big(G^{H_3}_{0,n-1}\big)-\frac{\gamma^2-\gamma^{3/2}}{4\sqrt{5}}\big(\gamma+\sqrt{\gamma}\big)^n+\frac12\Big|
\leq\frac14-\frac{1}{20}\sqrt{5}\approx 0.1382.$$

In his PhD thesis \cite{Haggstrom1994}, H\"{a}ggstr\"{o}m gave a description of the uniform spanning tree on the simple ladder graph as a Markov chain. A similar description will be given for the helix-3-graph in Section~\ref{S6}. The Markov chain description allows for very efficient computer simulations of the random spanning tree. However, the rigorous description of the Markov chain is a bit technical and is therefore deferred to Section~\ref{S6}.
\subsubsection{The enhanced helix-3-graph $G^{H_3'}$}
\label{S1.3.4}
The uniform spanning tree on the helix-3-graph results in a renewal process of order 2. It is tempting to conjecture that - in the spirit of Theorems~\ref{T1} and \ref{T2} - every renewal process of order 2 that is a Hermitian determinantal process could be realized via a spanning tree on a helix-3-graph by assigning an edge weight $c>0$ to the edges $(z_i)$. However, it is clear that renewal processes of order 2 whose Fourier transforms are the inverses of trigonometric polynomials of order 2 have one more parameter than those of order one. Hence, we will need an enhancement of the helix-3-graph, the graph $G^{H_3'}$. This graph consists of the vertices and edges of $G_{H_3}$ but in addition there are edges $g_i$, $i\in\Z$; that connect the vertices $i-2$ and $i$.

We assign edge weights 1 to the edges $(h_i)$, weights $c>0$ to the edges $(z_i)$ and $d\geq0$ to the edges $(g_i)$.
\begin{figure}[h]
\begin{picture}(400,140)(10,00)
\put(20,02){\includegraphics[scale=1.2]{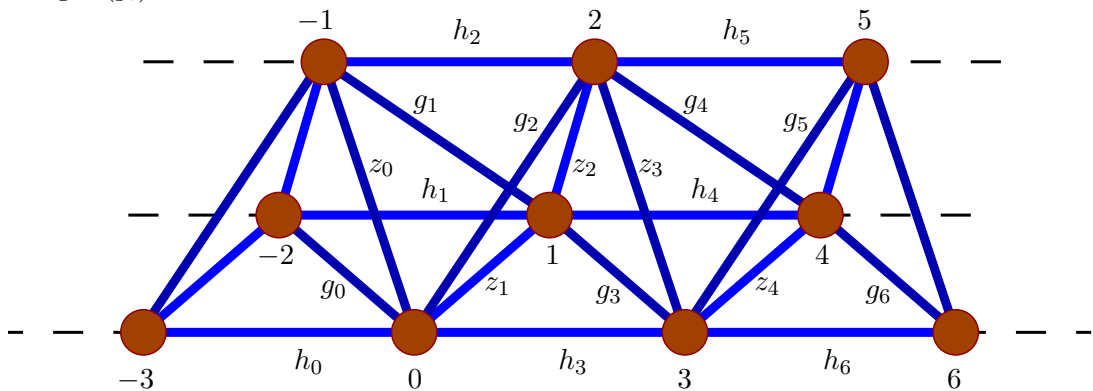}}
\put(063,06){$-3$}
\put(173,06){$0$}
\put(275,06){$3$}
\put(377,06){$6$}
\put(225,53){$1$}
\put(327,53){$4$}
\put(116,53){$-2$}
\put(131,142){$-1$}
\put(241,142){$2$}
\put(343,142){$5$}
\put(330,13){$h_6$}
\put(130,13){$h_{0}$}
\put(230,13){$h_{3}$}
\put(140,42){$g_{0}$}
\put(175,112){$g_{1}$}
\put(213,105){$g_{2}$}
\put(315,105){$g_{5}$}
\put(244,40){$g_{3}$}
\put(277,112){$g_{4}$}
\put(346,40){$g_{6}$}
\put(178,78){$h_{1}$}
\put(280,78){$h_{4}$}
\put(190,138){$h_{2}$}
\put(292,138){$h_{5}$}
\put(202,43){$z_{1}$}
\put(304,43){$z_{4}$}
\put(260,88){$z_{3}$}
\put(158,88){$z_{0}$}
\put(235,88){$z_{2}$}
\end{picture}
\caption[b]{A finite section of the enhanced helix-3-graph}
\label{F1.4}
\end{figure}

While $X_m=\1_T(z_m)$ clearly defines a renewal process of order 2, it is not the case that every such renewal process can be obtained by choosing $c$ and $d$ and a probability $p$ for a Bernoulli thinning appropriately. In fact, let $(Y_m)$ be a nontrivial (i.e., non-Bernoulli) renewal process of order one that is a determinantal process. Let $(Y'_m)$ be an independent copy and define $Z_{2m}=Y_m$ and $Z_{2m+1}=Y'_m$. Then $(Z_m)$ is a renewal process of order 2 and is determinantal. Now assume that there exist parameters $c,d$ and $p$ such that $(X_m\cdot B_m)$ equals $(Z_m)$ in distribution. Here $(B_m)$ is an independent family of Bernoulli random variables with parameter $p$. Since $Z_0$ and $Z_1$ are independent, also $X_0$ and $X_1$ are independent. This however implies, for the induced current that $I_{z_0,z_1}=0$. Recall that the current is induced by hooking a battery at $-1$ and $0$ such that the induced current along $z_0$ is 1. Now $I_{z_0,z_1}=0$ implies that the voltage at $1$ equals the voltage at $0$. However, since $c+d>0$ (otherwise the spanning tree is trivial), the voltage at any point $i\not\in\{0,1\}$ is strictly between the voltages at $-1$ and at $0$. This yields a contradiction.

\begin{theorem}
\label{T4}
\begin{enumerate}[(i)]
\item
For the random spanning tree on the enhanced helix-3-graph $G_{H_3'}$ with edge weights $c,d\geq0$, the process $X$ is a determinantal point process with (real-valued symmetric Toeplitz) matrix
$A_{k,l}=\hatf(l-k)$ given by
\begin{equation}
\label{E1.23}
\hatf(m)=\frac{x_1}{g'(x_1)}x_1^{|m|}+\frac{x_2}{g'(x_2)}x_2^{|m|},
\end{equation}
where
\begin{equation}
\label{E1.24}
g(x)=c^{-1}\big[x^4+(d+2)x^3+(c+2d+3)x^2+(d+2)x+1]
\end{equation}
and
$$x_{1}=\frac{1}{4}\big(-d-2+\sqrt {{d}^{2}-4\,d-4\,c}\big)+\frac14\sqrt { \left( -d
-2+\sqrt {{d}^{2}-4d-4\,c} \right) ^{2}-16}
$$
and
$$x_{2}=\frac{1}{4}\big(-d-2-\sqrt {{d}^{2}-4\,d-4\,c}\big)+\frac14\sqrt { \left( -d
-2-\sqrt {{d}^{2}-4d-4\,c} \right) ^{2}-16}
$$
are the two roots of $g$ in $\{z\in\C:\,|z|<1\}$.

 The Fourier transform $f$ of $\hatf$ equals
\begin{equation}
\label{E1.25}
f(x)=\frac{c}{(c+2d+3)+2(d+2)\cos(2\pi x)+2\cos(4\pi x)}.
\end{equation}
\end{enumerate}
\end{theorem}

Although the formula for $A_{k,l}$ is similar to that in Theorem~\ref{T3}, here $x_1$ and $x_2$ need not be complex conjugates. In fact, $x_1=\overline x_2$ if and only if $c\geq \frac{d^2}{4}-d$.

Note that the case $d+2=0$ corresponds to the counterexample with two interlaced and independent renewal processes. Since $d$ cannot assume negative values, this a more formal indication why not every renewal process of order 2 can be realized via a spanning tree.
\subsection{Outline}
\label{S1.4}
In Section~\ref{S2}, we give the proofs of Theorem~\ref{T1} and Proposition~\ref{P1.3}. We will present various approaches to the spanning tree on the ladder graph: an elementary and a slightly more sophisticated counting approach and the electrical network approach. In Section~\ref{S3}, we use the versatile counting approach from Section~\ref{S2} to prove Theorem~\ref{T2} and Proposition~\ref{P1.4}. In Section~\ref{S4}, we prove Theorem~\ref{T3}. Finally,  in Section~\ref{S6}, we present a Markov chain approach (Theorem~\ref{T5}) to the uniform spanning tree on the helix-3-graph in the spirit of H\"{a}ggstr\"{o}m \cite{Haggstrom1994}.

\Section{The Simple Ladder Graph}
\label{S2}
Recall the definition of the simple ladder graph $G^L$ and the finite subgraphs $G^L_{m,n}$ from Section~\ref{S1.3.1}. Define
$$V^{L,-}_{n}:=\{0,1\}\times\{n,n-1,\ldots\}\mbu V^{L,+}_m:=\{0,1\}\times\{m,m+1,\ldots\}$$
and the induced edge sets and subgraphs $E^{L,-}_{n}$, $G^{L,-}_{n}=(V^{L,-}_{n},E^{L,-}_n)$, $E^{L,+}_{m}$, and $G^{L,+}_{m}=(V^{L,+}_{m},E^{L,+}_m)$.

Recall that the weight function is defined on edges by $\weight(h_{i,m})=1$ and $\weight(z_m)=c$ for all $i=0,1$ and $m\in\Z$. Write $\sigma^c_n:=\weight\big(\SPT\big(G^L_{0,n-1}\big)\big)$ and note that for the weighted spanning tree measure $\P^c_n$ on $G^L_{-n,n}$, we have
\begin{equation}
\label{E2.01}
\P^c_n[\{t\}]=\frac{\displaystyle c^{\#\{k\in\{-n,\ldots,n\}:\,z_k\in t\}}}{\sigma^c_{2n+1}}.
\end{equation}

\subsection{The elementary counting approach}
\label{S2.1}
In this section, we present an elementary proof of Proposition~\ref{P1.3} and Theorem~\ref{T1}.\medskip

\textbf{Proof of Proposition~\ref{P1.3}.}\par
It is clear from \equ{E2.01} that we have to know the weighted number $\sigma^c_n$ of  spanning trees on $G^L_{0,n-1}$. To this end, we will derive a recursion formula.

Each spanning tree $t\in\SPT_{n+1}:=\SPT\big(G^L_{0,n}\big)$ falls into exactly one of the following four classes:
\begin{equation}
\label{defSTi}
\begin{aligned}
\SPT^1_{n+1}&:=\big\{t\in\SPT_{n+1}:\,h_{0,n},h_{1,n}\in t,\,z_n\not\in t\big\}\\
\SPT^2_{n+1}&:=\big\{t\in\SPT_{n+1}:\,h_{0,n},z_n\in t,\,h_{1,n}\not\in t\big\}\\
\SPT^3_{n+1}&:=\big\{t\in\SPT_{n+1}:\,h_{1,n},z_n\in t,\,h_{0,n}\not\in t\big\}\\
\SPT^4_{n+1}&:=\big\{t\in\SPT_{n+1}:\,h_{0,n},h_{1,n},z_n\in t\big\}.
\end{aligned}
\end{equation}
For $i=1,2,3$, the map $\SPT^i_{n+1}\to\SPT_n$, $t\mapsto t^*:=t\cap E^L_{0,n-1}$ is a bijection and
$$\weight(t)=\cases{\weight(t^*),&\mfalls i=1,\\[2mm]c\cdot\weight(t^*),&\mfalls i=2,3.}$$

For $t\in\SPT^4_{n+1}$, we have $z_{n-1}\not\in t$ and $t\cap E^L_{0,n-1}$ is not a tree but a forest that consists of two components, one linked to $(0,n-1)$ and the other linked to $(1,n-1)$. Hence
$$t^\dagger:=(t\cap E^L_{0,n-1})\cup\{z_{n-1}\}\in\SPT_n\setminus\SPT^1_{n}$$
and $\SPT^4_{n+1}\to\SPT_n\setminus\SPT^1_n$, $t\mapsto t^\dagger$ is a bijection. Clearly, we have $\weight(t)=\weight(t^\dagger)$ yielding
\begin{equation}
\label{E2.02}
\weight(\SPT^i_{n+1})=\cases{
\weight(\SPT_n),&\mfalls i=1,\\[2mm]
c\cdot \weight(\SPT_n),&\mfalls i=2,3,\\[2mm]
\weight(\SPT_n)-\weight(\SPT_{n-1}),&\mfalls i=4.}
\end{equation}

Recall $\sigma^c_n=\weight(\SPT_n)$. Hence \equ{E2.02} yields the recursion
\begin{equation}
\label{E2.03}
\sigma^c_{n+1}=(2c+2)\sigma^c_n-\sigma^c_{n-1}.
\end{equation}
Clearly, we have $\sigma^c_0=0$ and $\sigma^c_1=c$. In particular, for $c=1$, we have
\begin{equation}
\label{E2.04}
\sigma_0=0,\;\sigma_1=1\mbu\sigma_{n+1}=4\sigma_n-\sigma_{n-1}\mf n\in\N.
\end{equation}

 Hence we can compute
\begin{equation}
\label{E2.05}
\begin{array}{c|l|r}
n&\qquad\qquad\sigma^c_n&\sigma^1_n\\\hline
1&c&1\\
2&2c+2c^2&4\\
3&3c+8c^2+4c^3&15\\
4&4c+20c^2+24c^3+8c^4&56\\
5&5\,c+40\,{c}^{2}+84\,{c}^{3}+64\,{c}^{4}+16\,{c}^{5}&209\\
6&6\,c+70\,{c}^{2}+224\,{c}^{3}+288\,{c}^{4}+160\,{c}^{5}+32\,{c}^{6}&780
\end{array}
\end{equation}

The standard approach to an explicit solution of the recurrence relation \equ{E2.03} is to make the ansatz
$$\sigma^c_n=b_1\lambda_1^n+b_2\lambda_2^n$$
for certain $b_1,b_2,\lambda_1,\lambda_2\in\R$. From \equ{E2.03}, we get that $\lambda_{1,2}$ are the solutions of
$$\lambda^2=(2c+2)\lambda-1,$$
that is,
\begin{equation}
\label{E2.06}
\lambda_1:=c+1+\sqrt{c^2+2c}\mbu \lambda_2:=c+1-\sqrt{c^2+2c}.
\end{equation}
Note that
\begin{equation}
\label{E2.07}
\lambda_1=\frac{1}{\alpha}\mbu \lambda_2=\alpha,
\end{equation}
where $\alpha=c+1-\sqrt{c^2+2c}$ by \equ{E1.15}.
Comparing with the values $\sigma^c_1=c$ and $\sigma^c_2=2c+2c^2$ gives the coefficients
$$b_1=\frac12\frac{1-\alpha}{1+\alpha}\mbu b_2=-\frac12\frac{1-\alpha}{1+\alpha}.$$
Summing up, we have
$$\sigma^c_n=\frac12\,\frac{1-\alpha}{1+\alpha}\;\big(\alpha^{-n}-\alpha^n\big)
$$
which shows Proposition~\ref{P1.3}.
\eop

With the weighted number of spanning trees at hand, we are ready to prove Theorem~\ref{T1}.\medskip

\textbf{Proof of Theorem~\ref{T1}.}\par
Recall that $T$ is the weighted random tree under $\P^c$. Now we can compute, e.g.,
$$\P^c[z_0\in T]=\lim_{n\to\infty}\P^c_n[z_0\in T].$$
To this end, we deduce from \equ{defSTi} and \equ{E2.02} that
$$
\begin{aligned}
\weight\big(\big\{t\in\SPT\big(G^L_{-n,0}\big):\,z_0\in t\big\}\big)
&=\weight\big(\big\{t\in\SPT\big(G^L_{0,n}\big):\,z_0\in t\big\}\big)\\
&=\weight\big(\SPT_{n+1}\setminus\SPT^1_{n+1}\big)=\sigma^c_{n+1}-\sigma^c_{n}.
\end{aligned}$$
Hence
\begin{equation}
\label{E2.08}
\begin{aligned}
\P^c[z_0\in T]
&=\lim_{n\to\infty}\frac{\weight\big(\big\{t\in\SPT(G^{L}_{-n,n}):\,z_0\in t\big\}\big)}{\weight\big(\SPT(G^L_{-n,n})\big)}\\
=&\lim_{n\to\infty}\frac{\displaystyle c^{-1}\weight\big(\big\{t\in\SPT(G^{L}_{-n,0}):\,z_0\in t\big\}\big)\cdot\weight\big(\big\{t\in\SPT(G^{L}_{0,n}):\,z_0\in t\big\}\big)}{\weight\big(\SPT(G^L_{-n,n})\big)}\\
=&\lim_{n\to\infty}\frac{c^{-1}(\sigma^c_{n+1}-\sigma^c_{n})^2}{\sigma^c_{2n+1}}
=\frac{1-\alpha}{1+\alpha}\frac{(1-\alpha)^2}{2c\,\alpha}
=\frac{1-\alpha}{1+\alpha}.
\end{aligned}\end{equation}

In the same spirit, for $m\in \N$, we can compute $\P^c[z_0,z_m\in T]$. Let $$\SPT^+_n:=\{t\in\SPT_n:\,z_{n-1}\in t\}\mbu \SPT^-_n:=\{t\in\SPT_n:\,z_{0}\in t\}.$$
Note that $\SPT_n\setminus\SPT^+_n=\SPT^1_n$. Hence, by \equ{E2.02} and symmetry, we have
$$\weight(\SPT_n\setminus\SPT^-_n)=\weight(\SPT_n\setminus\SPT^+_n)=\sigma^c_{n-1}.$$
Arguing as in the derivation of \equ{E2.02}, we get
$$\weight\big(\SPT_n\setminus(\SPT^-_n\cup\SPT^+_n)\big)=\sigma^c_{n-2}.$$
Hence the inclusion/exclusion principle together with \equ{E2.03} yield
$$\sigma'_n:=\weight\big(\SPT^-_n\cap\SPT^+_n\big)=\sigma^c_n-2\sigma^c_{n-1}+\sigma^c_{n-2}=2c\,\sigma^c_{n-1}.$$
Similarly as in \equ{E2.08}, we conclude
\begin{equation}
\label{E2.09}
\begin{aligned}
\P^c[z_0&,z_m\in T]
=\lim_{n\to\infty}\frac{\weight\big(\big\{t\in\SPT(G^L_{-n,n}):\,z_0\in t\big\}\big)}{\weight(\SPT(G^L_{-n,n}))}\\
&=\lim_{n\to\infty}\frac{c^{-1}\weight\big(\big\{t\in\SPT(G^L_{-n,0}):\,z_0\in t\big\}\big)\weight\big(\big\{t\in\SPT(G^L_{0,n}):\,z_0\in t\big\}\big)}{\weight(\SPT(G^L_{-n,n}))}\\
&=\lim_{n\to\infty}\frac{c^{-2}(\sigma^c_{n+1}-\sigma^c_{n})\,2c\,\sigma^c_{m}\,(\sigma^c_{n+1-m}-\sigma^c_{n-m})}{\sigma^c_{2n+1}}\\
&=\left(\frac{1-\alpha}{1+\alpha}\right)^2\,\big(1-\alpha^{2m}\big).
\end{aligned}
\end{equation}
Iterating the argument, for $k=2,3,\ldots$ and $m_1<m_2<\ldots<m_k$, we get
\begin{equation}
\label{E2.10}
\P^c[z_{m_1},\ldots,z_{m_k}\in T]=\left(\frac{1-\alpha}{1+\alpha}\right)^k\,\prod_{l=2}^k\Big(1-\alpha^{2(m_l-m_{l-1})}\Big).
\end{equation}
By induction, we get that this probability is the determinant of the matrix $$\Big(\frac{1-\alpha}{1+\alpha}\;\alpha^{|m_{i}-m_{j}|}\Big)_{i,j=1,\ldots,k}.$$
Hence, $(X_m)_{m\in\Z}$ is in fact a determinantal point process with matrix $A$ given by
\equ{E1.12}.

The renewal property of $X$ follows from \cite[Proposition 2.10]{LyonsSteif2003}, but it can also be derived in the context of the spanning trees in a simple and intuitive way:
There are exactly $2m$ spanning trees $t\in\SPT\big(G^L_{\{0,\ldots,m\}}\big)$ with $z_0,z_m\in t$ but $z_k\not\in t$ for all $k=1,\ldots,m-1$. In fact, each of these trees contains all edges of the type $h_{i,k}$, $k=1,\ldots,m$, $i=0,1$, but one. Each of these trees has weight $c^2$.
Hence (recall $r_m$ from \equ{E1.08})
\begin{equation}
\label{E2.11}
\begin{aligned}
\P^c\big[T\cap\{z_0,&\ldots,z_m\}=\{z_0,z_m\}\big]\\
&=\lim_{n\to\infty}(\sigma^c_{2n+1}c^{2})^{-1}\weight\big(\big\{t\in\SPT(G^L_{-n,0}):\,z_0\in t\big\}\big)\\
&\qquad\qquad\times
\weight\big(\big\{t\in\SPT_{m+1}:\,t\cap\{z_0,\ldots,z_m\}=\{z_0,z_m\}\big\}\big)\\
&\qquad\qquad\times
\weight\big(\big\{t\in\SPT(G^L_{m,n}):\,z_m\in t\big\}\big)\\[2mm]
&=\lim_{n\to\infty}
\frac{c^{-2}(\sigma^c_{n+1}-\sigma^c_{n})\,2m\,c^2\,(\sigma^c_{n+1-m}-\sigma^c_{n-m})}{\sigma^c_{2n+1}}\\
&=\frac{(1-\alpha)^3}{1+\alpha\phantom{)}}\,m\,\alpha^{m-1}=r_m\,\P^c[z_0\in T].
\end{aligned}
\end{equation}
This finishes the proof of Theorem~\ref{T1}.\eop

\subsection{The more systematic counting approach}
\label{S2.2}
While the counting approach in Section~\ref{S2.1} worked out well, with a view to more complicated graphs, we present a more versatile method for counting the weighted spanning trees.

The main task is to count the number $\sigma^c_n$ of spanning trees on $G^L_{0,n-1}$. It became clear in Section~\ref{S2.1} that there is no linear recursion for $\sigma^c_n$ of first order, but rather of second order (see \equ{E2.03}). We aim at a more systematic derivation of the recursion formula \equ{E2.03} and its solution \equ{E1.16}.

For $m\leq n$, define $\Gamma^-_{m,n;1}:=\SPT\big(G^L_{m,n}\big)$ as the set of spanning trees on $G^L_{m,n}$ and let
 $$\Gamma^-_{m,n;2}:=\big\{t\subset E^L_{m,n}\setminus\{z_n\}:\,t\cup\{z_n\}\in\SPT\big(G^L_{m,n}\big)\big\}.$$
That is, the elements of $\Gamma^-_{m,n;2}$ are spanning forests with two connected components each of which contains exactly one of the points $(0,n)$ and $(1,n)$. Furthermore, let $\Gamma^+_{m,n;i}$ denote the image of $\Gamma^-_{m,n;i}$ under the map $k\mapsto m+n-k$. That is, $\Gamma^+_{m,n;1}=\Gamma^-_{m,n;1}$ and
$$\Gamma^+_{m,n;2}:=\big\{t\subset E^L_{m,n}\setminus\{z_m\}:\,t\cup\{z_m\}\in\SPT\big(G^L_{m,n}\big)\big\}.$$

Similarly, define $\Gamma^-_{n;i}$ and $\Gamma^+_{m;i}$, $i=1,2$, for the one-sided infinite graphs $G^{L,-}_n$ and $G^{L,+}_m$. Furthermore, we define
$$\Gamma^\pm_{n,n-1;1}=\emptyset \mbu \Gamma^\pm_{n,n-1;2}=\{\emptyset\}.$$
Define
$$a_{n;i}=a^c_{n;i}:= \weight\big(\Gamma^-_{0,n-1;i}\big)$$
and $a_n:=(a_{n;1},a_{n;2})$. Note that $a_0=(0,1)$.

For $i,j=1,2$ and $t\in\Gamma^-_{0;i}$, let
\begin{equation}
\label{E2.12}
F_{i,j}:=\big\{F\subset\{h_{0,1},\,h_{1,1},\,z_{1}\}:\,t\cup F\in\Gamma^-_{1;j}\big\}
\end{equation}
and define
\begin{equation}
\label{E2.13}
M_{i,j}:=\weight\big(F_{i,j}\big).
\end{equation}
That is, we have
$$\begin{aligned}
F_{1,1}&=\big\{\{h_{0,1},h_{1,1}\},\,\{h_{0,1},z_1\},\,\{h_{1,1},z_1\}\big\},\\
F_{1,2}=F_{2,1}&=\big\{\{h_{0,1}\},\,,\{h_{1,1}\}\big\},\\
F_{2,2}&=\big\{\{h_{0,1},h_{1,1}\}\big\}.
\end{aligned}
$$
Note that $F_{1,1}$ corresponds to the cases $\SPT^i_{n+1}$, $i=1,2,3$ and $F_{2,1}$ corresponds to $\SPT^4_{n+1}$.

Going through all possibilities for $F\in F_{i,j}$, we get the matrix
\begin{equation}
\label{E2.14}
M=\left(
\begin{array}{cc}
2c+1& 2\\
c& 1
\end{array}
\right).
\end{equation}
Similarly, let
$F'_{i,j}:=\big\{F\in F_{i,j}:\,z_{1}\in F\big\}$
and
$M'_{i,j}:=\weight\big(F'_{i,j}\big)$.
Then
\begin{equation}
\label{E2.15}
M'=\left(
\begin{array}{cc}
2c& 0\\
c& 0
\end{array}
\right)
\end{equation}

Now instead of \equ{E2.03}, the recursion reads
$a_0=(0,1)$ and
\begin{equation}
\label{E2.16}
a_{n+1}=a_n\,M\mf n\geq0.
\end{equation}
Hence
\begin{equation}
\label{E2.17}
a_n=(0,1)\cdot M^{n}\mf n\in\N.
\end{equation}
The computation in \equ{E2.17} can be facilitated by an eigenvector decomposition. The characteristic polynomial of $M$ is
\begin{equation}
\label{E2.18}
\chi_M(x)=x^2-2(c+1)x+1.
\end{equation}
The roots are (recall $\alpha$ from \equ{E1.15})
$\lambda_1=1/\alpha$ and $\lambda_2=\alpha$.
The corresponding left eigenvectors $w_1$ and $w_2$ are the row vectors of the matrix $W=(w_{ij})_{i,j=1..2}$ given by
\begin{equation}
\label{E2.19}
W=\left(\begin{array}{cc}
(\alpha^{-1}-1)/2&1\\
(\alpha-1)/2&1
\end{array}\right)
\end{equation}
Let
\begin{equation}\label{E2.20}
C:=a_0\cdot W^{-1}
=\left(\frac{\alpha}{1+\alpha},\,\frac{1}{1+\alpha}\right).
\end{equation}
Summarizing the discussion, we get \begin{equation}
\label{E2.21}
\sigma^c_n=a_{n;1}
=\sum_{j=1}^2C_j\,w_{j,1}\,\lambda_j^{n}\;=\;\frac12\frac{1-\alpha}{1+\alpha}\big(\alpha^{-n}-\alpha^{n}\big).
\end{equation}
Hence we have shown Proposition~\ref{P1.3} by means of a more systematic counting approach than in Section~\ref{S2.1}.

We now come to a more systematic approach to probabilities such as in \equ{E2.08}.
For $m\in\{-1,-2,\ldots\}\cup\{-\infty\}$, $n\in\{2,3,\ldots\}\cup\{\infty\}$, $i,j=1,2$,  $t^-\in\Gamma^-_{m,0;i}$, and $t^+\in\Gamma^+_{1,n;j}$ define the sets
$$B_{i,j}:=\Big\{B\subset \big\{h_{0,1},\,h_{1,1}\big\}:\,t^-\cup B\cup t^+\in\SPT\big(G^L_{m,n}\big)\Big\}.
$$
Note that the set $B_{i,j}$ of ``bridges'' between $t^-$ and $t^+$ depends only on the values of $i$ and $j$ but not on $m,n$ and the choices of $t^-$ and $t^+$. Now define the matrix of weights
$$N_{i,j}=\weight(B_{i,j})=\#B_{i,j}.$$
Clearly
$$
B_{1,1}=\big\{\{h_{0,1}\},\,\{h_{1,1}\}\big\},\quad
B_{1,2}=B_{2,1}=\big\{\{h_{0,1},h_{1,1}\}\big\}\mbu
B_{2,2}=\emptyset.
$$

Simple counting gives the values
\begin{equation}
\label{E2.22}
N:=\left(\begin{array}{cc}
2&1\\
1&0
\end{array}\right).
\end{equation}
Note that
$$\weight\big(\{t\in\Gamma^-_{-n,0;1}:\,z_0\in t\}\big)=(a_{n}M')_1$$
and that $(a_{n}M')_2=0$.
Hence
\begin{equation}
\label{E2.23}
\P^c_n[z_0\in T]=\frac{a_{n}\,M'N\,(a_{n})^T}{\weight(\SPT(G^L_{-n,n}))}
=\frac{a_{n}M'N(a_{n})^T}{a_{n}\,M\,N\,(a_{n})^T}.
\end{equation}
Letting $n\to\infty$, only the eigenvector $w_1$ for the largest eigenvalue $\lambda_1$ in the representation \equ{E2.21} contributes. Hence
\begin{equation}
\label{E2.24}
\P^c[z_0\in T]=\lim_{n\to\infty}\P^c_n[z_0\in T]=\frac{w_1M'N(w_1)^T}{w_1MN(w_1)^T}
=\frac{w_1M'N(w_1)^T}{\lambda_1\,w_1N(w_1)^T}.
\end{equation}
With the explicit values from \equ{E2.06}, \equ{E2.19} and \equ{E2.22}, we obtain
$$\P^c[z_0 \in T]=\frac{1-\alpha}{1+\alpha}.$$
Iterating the argument, we also get (for $m\in \N$)
$$\P^c_n[z_0,z_m\in T]=\frac{a_{n}M'M^{m-1}N(M')^Ta_{n-m}^T}{a_{n}M^{m+1}Na_{n-m}^T}.$$
Letting $n\to\infty$, we get
\begin{equation}
\label{E2.25}
\P^c[z_0,z_m\in T]=\frac{w_1M'M^{m-1}N(M')^T(w_1)^T}{\lambda_1^{m+1}\,w_1N(w_1)^T}.
\end{equation}
Let $(\mu_1,\mu_2)=w_1\,M'\,W^{-1}$
be the coefficients of the eigenvector (for $M$) decomposition of $w_1\,M'$; that is, $w_1\,M'=\mu_1w_1+\mu_2w_2$. Then
$$\P^c[z_0,\,z_m\in T]=\frac{\sum_{j=1}^2\mu_j\,\lambda_j^{m-1}\,w_j\,N\,(w_1M')^T}{\lambda_1^{m+1}\,w_1\,N\,(w_1)^T}.$$
We define
$\tilde\lambda_j=\lambda_j/\lambda_1$, that is $\tilde\lambda_1=1$ and $\tilde\lambda_2=\alpha^2$, and
$$\varrho_j:=\frac{\mu_j\,w_j\,N\,(w_1M')^T}{\lambda_1\lambda_j\,w_1\,N\,(w_1)^T}\mf j=1,2,$$
that is, $\varrho_1=\big((1-\alpha)/(1+\alpha)\big)^2$ and $\varrho_2=-\varrho_1$.
Then
\begin{equation}
\label{E2.26}
\P^c[z_0,z_m\in T]=\sum_{j=1}^2\varrho_j\,\tilde\lambda_j^m\,=\,\left(\frac{1-\alpha}{1+\alpha}\right)^2\big(1-\alpha^{2m}\big).
\end{equation}
Since $\P[z_0,z_m\in T]=(A_{0,0})^2-(A_{0,m})^2$, from \equ{E2.26} we get the values $|A_{k,l}|$ for all $k,l\in\Z$. In order to infer the signs of the matrix entries, we could compute the probabilities for three edges to be in the tree and then proceed as for the helix-3-graph (see Section~\ref{S4.2}). For the simple ladder graph, however, the electric network approach that we present in the next section shows in a simple way that all entries of $A$ can be chosen to be nonnegative.

Concluding, we have developed a matrix-based counting procedure that enables us to compute probabilities for the random spanning tree $T$. While for the case of the simple ladder graph this method might seem to be a bit exaggerated, it simplifies things when we come to different graphs and essentially have to compute the matrices $M$, $M'$ and $N$ as well as the eigenvalues and left eigenvectors of $M$.
\subsection{The electrical network approach}
\label{S2.3}
The aim of this section is to give a proof of Theorem~\ref{T1} that, in the spirit of \cite{BurtonPemantle1993}, relies on the interpretation of the matrix $A$ as the transfer current matrix.

Consider the simple ladder graph $G^L$ with the weight function $\weight$ as an electrical network with conductance $\weight(e)$ along the edge $e$. We will need to give an orientation to the edges $z_m$ and we choose the orientation from $(1,m)$ to $(0,m)$.
Recall that $G^{L,-}_m$ and $G^{L,+}_m$ are the subgraphs of $G^L$ with vertex sets
$$V^{L,-}_{m}:=\{0,1\}\times\{m,m-1,\ldots\}\mbu V^{L,+}_m:=\{0,1\}\times\{m,m+1,\ldots\},$$
respectively.

 Let $r$ and $r^+$ denote the effective resistances between the sites $(0,0)$ and $(1,0)$ in the networks $G^L$ and $G^{L,+}_0$, respectively. Furthermore, let $r_0$ denote the effective resistance between $(0,0)$ and $(1,0)$ in $G^{L,+}_0\setminus\{z_0\}$. By serial connection, we get $r_0=r^++2$. Considering the two parallel resistors $z_0$ and $G^{L,+}_0\setminus\{z_0\}$, we get
$$r^+=\frac{1}{\displaystyle c+\frac{1}{r_0}}=\frac{1}{\displaystyle c+\frac{1}{r^++2}}.$$
The solution is $r^+=(1-\alpha)/c$.
Considering the three parallel resistors $G^{L,-}_0\setminus\{z_0\}$, $G^{L,+}_0\setminus\{z_0\}$ and $z_0$, we get
\begin{equation}
\label{E2.27}
r=\frac{1}{c+\frac{2}{r_0}}=c^{-1}\,\frac{1-\alpha}{1+\alpha}.
\end{equation}
By Kirchhoff's theorem \equ{E1.01}, this yields $\P[z_0\in T]=r\cdot\weight(z_0)=\frac{1-\alpha}{1+\alpha}$.

The next task is to compute the full transfer current matrix. Connect a battery at $(0,0)$ and $(1,0)$ with voltages $u(0,0)=0$ and $u(1,0)=r$. The resulting electric current is a unit flow through the network. Denote by $u(i,k)$ the voltage at $(i,k)$ and let $u(k):=u(1,k)-u(0,k)$.

Let $l\geq0$. If we know $u(l)$, then we get $u(l+1)$ as follows: Consider the three serial resistors $h_{0,l+1}$, $h_{1,l+1}$ and $G^{L,+}_{l+1}$. Each of the resistors has a potential difference proportional to its resistance and the differences add up to $u(l)$. Hence
$$u(l+1)=\frac{r_+}{2+r_+}\,u(l)\;=\;\alpha\,u(l).$$

Together with the initial value $u(0)=r$, we get for $l\in\Z$,
$$u(l)=c^{-1}\,\frac{1-\alpha}{1+\alpha}\,\alpha^{|l|}.$$

By Ohm's law, the current along $z_l$ equals
$$I_{z_0,z_l}=u(l)\cdot\weight(z_l)=\frac{1-\alpha}{1+\alpha}\,\alpha^{|l|}.$$
By the Burton-Pemantle theorem, translation invariance and symmetry, we get $A_{k,l}=I_{0,|l-k|}$ and hence we have another proof for Theorem~\ref{T1}.

\Section{The zigzag ladder graph}
\label{S3}
Recall the definition of $G^{H_k}$ and $G^{H_k}_{m,n}$ from Section~\ref{S1.3.2}. Similarly as for the simple ladder graph, also define
$$V^{H_k,-}_n:=\big\{n,n-1,\ldots\big\}\mbu V^{H_k,+}_m:=\big\{m,m+1,\ldots\big\}$$
 and the induced subgraphs $G^{H_k,-}_{n}=\big(V^{H_k,-}_n, E^{H_k,-}_n\big)$  and $G^{H_k,+}_{m}=\big(V^{H_k,+}_m, E^{H_k,+}_m\big)$.

Here we focus on the case $k=2$.

\subsection{Proof of Proposition~\ref{P1.4}}
\label{S3.1}
For $m<n$, let $\Gamma^-_{m,n;1}=\SPT\big(G^{H_2}_{m,n}\big)$ denote the set of spanning trees on $G^{H_2}_{m,n}$ and let
$$\Gamma^-_{m,n;2}:=\Big\{t\subset E^{H_2}_{m,n}\setminus\{z_n\}:\,t\cup\{z_n\}\in\SPT\big(G^{H_2}_{m,n}\big)\Big\}.$$
That is, the elements of $\Gamma^-_{m,n;2}$ are spanning forests with two connected components each of which contains exactly one of the points $n-1$ and $n$. Similarly as in Section~\ref{S2.2} define the sets $\Gamma^\pm_{m;i}$ and $\Gamma^\pm_{m,n;i}$. Furthermore, define $\Gamma^\pm_{n,n;1}=\emptyset$ and $\Gamma^\pm_{n,n;2}=\{\emptyset\}$. Finally, let $a_{n;i}=\weight(\Gamma^-_{0,n-1;i})$ for $i=1,2$, $n\in\N$, and $a_n:=(a_{n;1},a_{n;2})$.

Recall that the weight function is defined by $\weight(h_m)=1$ and $\weight(z_m)=c$ for all $m\in\Z$ and that $\P^c$ is the weighted spanning tree distribution on $\SPT\big(G^{H_2}\big)$. We shall exhibit the flexibility of the method from Section~\ref{S2.2} by computing the probabilities $\P^c[z_0\in T]$ and $\P^c[z_0,z_m\in T]$.

Similarly as in Section~\ref{S2.2}, for $i,j=1,2$ and $t^-\in\Gamma^-_{0;i}$, $t^+\in\Gamma^+_{0;j}$ let
\begin{equation}
\label{E3.01}
F_{i,j}:=\big\{F\subset\{h_{1},z_{1}\}:\,t^-\cup F\in\Gamma^-_{1;j}\big\},
\end{equation}
\begin{equation}
\label{E3.02}
F'_{i,j}:=\big\{F\in F_{i,j}:\,z_1\in F\big\},
\end{equation}
and
\begin{equation}
\label{E3.03}
B_{i,j}:=\big\{B\subset \{h_1\}:\,t^-\cup B\cup t^+\in\SPT\big(G^{H_2}\big)\big\}.
\end{equation}
Note that independently of the actual choice of $t^-$ and $t^+$, we have
$$F_{1,1}=\big\{\{h_{1}\},\{z_{1}\}\big\},\quad
F_{1,2}=\{\emptyset\},\quad
F_{2,1}=\big\{\{h_{1},z_{1}\}\big\},\quad
F_{2,2}=\big\{\{h_{1}\}\big\}
$$
and
$$F'_{1,1}=\big\{\{z_{1}\}\big\},\quad
F'_{2,1}=\big\{\{h_{1},z_{1}\}\big\},\quad
F'_{1,2}=F'_{2,2}=\emptyset
$$
and
$$B_{1,1}=\{\emptyset\},\quad
B_{1,2}=B_{2,1}=\big\{\{h_{1}\}\big\},\quad
B_{2,2}=\emptyset.
$$

Define the matrices $M$, $M'$ and $N$ by
\begin{equation}
\label{E3.04}
M_{i,j}:=\weight(F_{i,j}),\qquad M'_{i,j}:=\weight(F'_{i,j})\mbu N_{i,j}=\weight(B_{i,j}).\end{equation}
Explicit counting yields
\begin{equation}
\label{E3.05}
M=\left(\begin{array}{cc}
1+c&1\\c&1
\end{array}\right),\qquad
M'=\left(\begin{array}{cc}
c&0\\c&0
\end{array}\right),\mbu
N=\left(\begin{array}{cc}
1&1\\1&0
\end{array}\right).
\end{equation}
Hence, we have the recursion $a_1=(0,1)$ and $a_{n+1}=a_nM$, that is
\begin{equation}
\label{E3.06}
a_n=(0,1)\cdot M^{n-1}\mf n\geq1.
\end{equation}

The characteristic polynomial of $M$ is
$$\chi_M(x)=x^2-(2+c)x+1$$
with roots (recall $\alpha$ from \equ{E1.17})
$\lambda_1=1/\alpha$ and $\lambda_2=\alpha$.
The corresponding left eigenvectors $w_1$ and $w_2$ are the row vectors of the matrix $W=(w_{ij})_{i,j=1..2}$ given by
\begin{equation}
\label{E3.07}
W=\left(\begin{array}{cc}
\alpha^{-1}-1&1\\
\alpha-1&1
\end{array}\right).
\end{equation}
Arguing as in \equ{E2.21}, with $C=a_1W^{-1}=\big(\alpha/(1+\alpha),1/(1+\alpha)\big)$, we get
\begin{equation}
\label{E3.08}
\sigma^c_n=a_{n;1}=\sum_{j=1}^2C_j\,w_{j,1}\,\lambda_j^{n-1}
\;=\;\frac{1-\alpha}{1+\alpha}\big(\alpha^{1-n}-\alpha^{n-1}\big).
\end{equation}
This finishes the proof of Proposition~\ref{P1.4}.\eop

\subsection{Proof of Theorem~\ref{T2}}
\label{S3.2}
Having the counting scheme for the zigzag ladder graph $G^{H_2}$ at hand, we can compute (recall \equ{E2.24})
\begin{equation}
\label{E3.09}
\P^c[z_0\in T]=\frac{w_1M'N(w_1)^T}{\lambda_1\,w_1N(w_1)^T}=\frac{1-\alpha}{1+\alpha}=\sqrt{\frac{c}{c+4}\,}.
\end{equation}
In particular, for the uniform spanning tree,
\begin{equation}
\label{E3.10}
\P^1[z_0\in T]=\frac{1}{\sqrt{5}}.
\end{equation}
Furthermore, for $m\in\N$, arguing as in \equ{E2.25}, we get
\begin{equation}
\label{E3.11}
\P^c[z_0,z_m\in T]=\frac{w_1M'M^{m-1}N(M')^T(w_1)^T}{\lambda_1^{m+1}\,w_1N(w_1)^T}.
\end{equation}
Arguing as in Section~\ref{S2.2}, we get
\begin{equation}
\label{E3.12}
\P^c[z_0,z_m\in T]=\left(\frac{1-\alpha}{1+\alpha}\right)^2\big(1-\alpha^{2m}\big).
\end{equation}
In particular, for the uniform spanning tree,
\begin{equation}
\label{E3.13}
\P^1[z_0,z_m\in T]=\frac15\big(1-\gamma^{-4m}\big),
\end{equation}
where $\gamma=(\sqrt{5}+1)/2$ is the golden ratio.

From \equ{E3.12}, we infer $|A_{k,l}|=\frac{1-\alpha}{1+\alpha}\,\alpha^{|l-k|}$
but we do not get the signs of the matrix entries. While $A_{k,k}>0$ is clear, the signs of the other entries need a little thought. We could compute the probabilities for three distinct edges to be in the spanning tree and infer the signs of the matrix (up to the free choice of the sign of $A_{0,1}$). We will do so for the helix-3-graph in Section~\ref{S4.2}. However, for the zigzag graph a different approach works out simpler. Using the Burton-Pemantle theorem of \cite{BurtonPemantle1993}, we interpret $A$ as the impedance matrix of an electrical network and compute the signs of the induced currents.

We start with the electrical network approach. Let us orient the edges $z_m$ from $m-1$ to $m$ and hook a battery at $-1$ with voltage $u(-1)=c^{-1}(1-\alpha)/(1+\alpha)$ and at $0$ with $u(0)=0$. By Kirchhoff's theorem and \equ{E3.09}, the resulting electric current is a unit flow from $-1$ to $0$. Denote by $I(k)=c(u(k-1)-u(k))$ the resulting current along $z_k$. We have $I(0)=(1-\alpha)/(1+\alpha)>0$. Clearly, $u(1)\geq0$ and hence $I(1)\leq0$. Continuing to the right, we see that
$$\min\big\{u(k-1),u(k)\big\}\,\leq\, u(k+1)\,\leq\, \max\big\{u(k-1),u(k)\big\}.$$
Hence $I(k+1)$ and $I(k)$ have different signs. Hence
$A_{0,k}=(-1)^k|A_{0,k}|=\frac{1-\alpha}{1+\alpha}(-\alpha)^{|k|}$ and thus
$$A_{k,l}=\frac{1-\alpha}{1+\alpha}\,(-\alpha)^{|l-k|}.$$

Note that $A$ and $A^{1/2}_{k,l}:=(-1)^{l-k}A_{k,l}$ define the same determinantal point process. Hence Theorem~\ref{T2} is proved.\eop

\ignore{As long as we consider events that involve only two edges, the signs of off-diagonal entries of $A$ do not matter. But as soon as we have three edges, the signs become relevant. In fact, if you do not like the the electrical network argument given above, you could derive the signs as follows (up to the arbitrary choice of the sign of $A_{0,1}$). Assume $m,n\geq1$ and compute
$$\P^c[z_0,z_m,z_{m+n}\in T]=\frac{w_1M'M^{m-1}M'M^{n-1}N(M')^T(w_1)^T}{\lambda_1^{m+n+1}\,w_1N(w_1)^T}.
$$
With an eigenvalue expansion as in Section~\ref{S2.2}, we can compute
$$\begin{aligned}
\P^c[z_0,z_m,z_{m+n}\in T]&=\left(\frac{1-\alpha}{1+\alpha}\right)^3\Big[
1-\alpha^{2m}-\alpha^{2n}+\alpha^{2m+2n}\Big]\\[2mm]
&=\left(\frac{1-\alpha}{1+\alpha}\right)^3\det\left(
\begin{array}{ccc}
1           &(-\alpha)^m& (-\alpha)^{m+n}\\
(-\alpha)^m    &       1& (-\alpha)^n    \\
(-\alpha)^{m+n}&(-\alpha)^n&        1
\end{array}\right).
\end{aligned}
$$
}
\begin{remark}\label{R3.1}
Using the transfer current idea, the signs of the matrix $A$ were easy to determine for the ladder graph and for the zigzag graph. For the helix-3-graph that will be studied in the following section, the electrical network is rather involved and the signs of $A$ do not follow (such) simple patterns. In fact, it is the sign of the real part of some rotation on the unit circle $S_1\subset\C$. See \equ{E1.21}.
\end{remark}
\Section{The helix-3-graph}
\label{S4}
In this section, we use the counting scheme from Section~\ref{S2.2} to prove Proposition~\ref{P1.5} and Theorem~\ref{T3}.

Recall the definitions of $G^{H_3}$, $G_{m,n}^{H_3}$, $G_{m}^{H_3,+}$ and $G_{n}^{H_3,-}$ from Section~\ref{S1.3.2} and Section~\ref{S3} and recall that $\weight\equiv1$.

\subsection{Counting the number of spanning trees: Proof of Proposition~\ref{P1.5}}
\label{S4.1}

Define $\Gamma^-_{m,n;1}=\SPT\big(G^{H_3}_{m,n}\big)$ as the set of spanning trees on $G^{H_3}_{m,n}$. For the simple ladder graph (Section~\ref{S2.2}) and for the zigzag graph (Section~\ref{S3.1}), in order to set up a first order recursion for the number of spanning trees, we introduced the sets of spanning forests that would become trees if we connected the two rightmost points. For the helix-3-graph, there are three rightmost points and we have to distinguish four types of spanning forests.

For $n\geq m+2$, we introduce the following subgraphs of $G^{H_3}_{m,n}$:
\begin{itemize}
\item
Let $\Gamma^-_{m,n;2}$ be the set of spanning forests with two connected components $C_1$ and $C_2$ such that $n\in C_1$ and $n-1,n-2\in C_2$. That is, by adding an extra bond between $n$ and either $n-1$ or $n-2$, we get a spanning tree.
\item
Define $\Gamma^-_{m,n;3}$ similarly as $\Gamma^-_{m,n;2}$ but with $n-2\in C_1$ and $n,n-1\in C_2$.
\item
Define $\Gamma^-_{m,n;4}$ similarly as $\Gamma^-_{m,n;2}$ but with $n-1\in C_1$ and $n,n-2\in C_2$.
\item
Let $\Gamma^-_{m,n;5}$ be the set of spanning forests with three connected components each of which contains exactly one of the points $n,n-1, n-2$.
\end{itemize}
More formally, we could describe these sets as follows.
$$\begin{aligned}
\Gamma^-_{m,n;2}&=\big\{t\subset E^{H_3}_{m,n}\setminus\{z_n\}:t\cup\{z_n\}\in\Gamma^-_{m,n;1}
  \mbs{and}t\cup\{z_{n+1},h_{n+1}\}\in\Gamma^-_{m,n+1;1}\big\},\\[1.5mm]
\Gamma^-_{m,n;3}&=\big\{t\subset E^{H_3}_{m,n}\setminus\{z_{n-1}\}:t\cup\{z_{n-1}\}\in\Gamma^-_{m,n;1}
  \mbs{and}t\cup\{z_{n+1},h_{n+1}\}\in\Gamma^-_{m,n+1;1}\big\},\\[1.5mm]
\Gamma^-_{m,n;4}&=\big\{t\subset E^{H_3}_{m,n}\setminus\{z_{n-1},z_n\}:t\cup\{z_{n-1}\}\in\Gamma^-_{m,n;1}
  \mbs{and}t\cup\{z_{n}\}\in\Gamma^-_{m,n;1}\big\},\\[1.5mm]
\Gamma^-_{m,n;5}&=\big\{t\subset E^{H_3}_{m,n}\setminus\{z_{n-1},z_n\}:t\cup\{z_{n-1},z_n\}\in\Gamma^-_{m,n;1}\big\}.
\end{aligned}
$$
For $m=-\infty$, we simply write $\Gamma^-_{n;i}$ for $\Gamma^-_{-\infty,n;i}$. Similarly, define $\Gamma^+_{m,n;i}$ as the image of $\Gamma^-_{m,n;i}$ of the map $\{m,\ldots,n\}\to\{m,\ldots,n\}$, $k\mapsto m+n-k$ that reverses the order of the vertices and define $\Gamma^+_{m;i}:=\Gamma^+_{m,\infty;i}$.

For $n=m$ and $n=m+1$, we pretend that the points $m-2$, $m-1$ and $m$ are disconnected and define
$$\Gamma^-_{m,m;1}:=\Gamma^-_{m,m;2}:=\Gamma^-_{m,m;3}:=\Gamma^-_{m,m;4}:=\emptyset\mbu
\Gamma^-_{m,m;5}:=\{\emptyset\}$$
as well as
$$\Gamma^-_{m,m+1;1}:=\Gamma^-_{m,m+1;2}:=\Gamma^-_{m,m+1;4}:=\emptyset,\quad
\Gamma^-_{m,m+1;3}:=\{z_{m+1}\},
\quad \Gamma^-_{m,m+1;5}:=\{\emptyset\}.$$
Define
$$a_{n;i}=\#\Gamma^-_{0,n-1;i},\qquad i=1,\ldots,5,$$
and $a_n:=(a_{n;1},a_{n;2},a_{n;3},a_{n;4},a_{n;5}).$
For $i,j=1,\ldots,5$, define $F$, $F'$, $M$, and $M'$ as in \equ{E3.04}.

\begin{figure}[t]
\label{F4.1}
\begin{picture}(400,120)(-50,0)
\put(20,02){\includegraphics[scale=0.9]{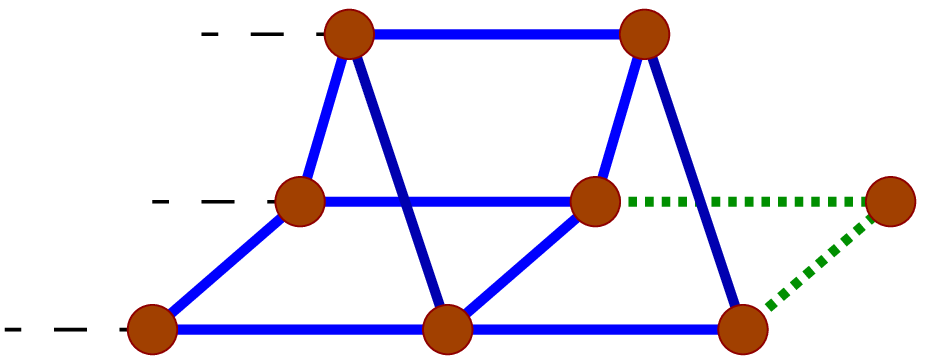}}
\put(126,00){$-3$}
\put(210,00){$0$}
\put(255,40){$1$}
\put(175,40){$-2$}
\put(180,110){$-1$}
\put(103,7){$h_{-3}$}
\put(172,7){$h_{0}$}
\put(140,60){$h_{-2}$}
\put(218,60){$h_{1}$}
\put(140,104){$h_{-1}$}
\put(158,30){$z_{-2}$}
\put(232,30){$z_{1}$}
\put(120,78){$z_{-3}$}
\put(199,78){$z_{0}$}
\put(161,78){$z_{-1}$}
\end{picture}
\caption[b]{Adding some of the edges $\{h_1,z_1\}$ changes how $0,-1,-2$ are connected.}
\end{figure}
\bigskip\par

\begin{table}[h]
$$
\def\arraycolsep{1.5mm}
\begin{array}{c|c|c|c|c|c|}
i\;\backslash\; j&1&2&3&4&5\\\hline
1&\{\{z_1\},\{h_1\}\}&\{\emptyset\}&\emptyset&\emptyset&\emptyset\\[1.5mm]
2&\{\{h_1,z_1\}\}&\emptyset&\{\{z_1\}\}&\{\{h_1\}\}&\{\emptyset\}\\[1.5mm]
3&\{\{h_1,z_1\}\}&\{\{h_1\}\}&\emptyset&\emptyset&\emptyset\\[1.5mm]
4&\emptyset&\emptyset&\{\{h_1\},\{z_1\}\}&\emptyset&\{\emptyset\}\\[1.5mm]
5&\emptyset&\emptyset&\{\{h_1,z_1\}\}&\emptyset&\{\{h_1\}\}\\\hline
\end{array}
$$
\caption{The sets $F_{i,j}$ for the helix-3-graph}
\label{Tab4.1}
\end{table}
Explicit counting (see Table~\ref{Tab4.1} for $F_{i,j}$) yields
\begin{equation}
\label{E4.01}
M=\left(
\begin{array}{ccccc}
2& 1& 0& 0& 0\\
1& 0& 1& 1& 1\\
1& 1& 0& 0& 0\\
0& 0& 2& 0& 1\\
0& 0& 1& 0& 1
\end{array}
\right)\mbu
M'=\left(
\begin{array}{ccccc}
1& 0& 0& 0& 0\\
1& 0& 1& 0& 0\\
1& 0& 0& 0& 0\\
0& 0& 1& 0& 0\\
0& 0& 1& 0& 0
\end{array}
\right)
.
\end{equation}
Now the recursion reads
$a_1=(0,0,0,0,1)$ and
\begin{equation}
\label{E4.02}
a_{n+1}=a_n\,M\mf n\geq0.
\end{equation}
Hence
\begin{equation}
\label{E4.03}
a_n=(0,0,0,0,1)\cdot M^{n-1}\mf n\in\N.
\end{equation}
 The characteristic polynomial of $M$ is
\begin{equation}
\label{E4.04}
\chi_M(x)=x^5-3x^4+3x-1=(x-1)\big(x^4-2x^3-2x^2-2x+1\big).
\end{equation}
Recall that $\gamma=(\sqrt{5}+1)/2$
denotes the golden ratio. The roots of $\chi_M$ are
\begin{equation}
\label{E4.05}
\begin{aligned}
\lambda_1&=\gamma+\sqrt{\gamma}\approx 2.890536\ldots,\\
\lambda_2&=1,\\
\lambda_3&=\gamma-\sqrt{\gamma}=1/\lambda_1\approx 0.346014\ldots,\\
\lambda_4&=-\frac1\gamma+\sqrt{-1/\gamma},\qquad (\mbsr{note that} |\lambda_4|=1),\\
\lambda_5&=-\frac1\gamma-\sqrt{-1/\gamma}=\bar\lambda_4.
\end{aligned}
\end{equation}
The corresponding left eigenvectors are the row vectors of the matrix $W=(w_{ij})_{i,j=1..5}$ given by
\begin{equation}
\label{E4.06}
W=\left(\begin{array}{ccccc}

\gamma^{5/2}+\gamma^2&\gamma+\sqrt{\gamma}&\gamma+1/\sqrt{\gamma}
&1&\gamma^{3/2}\\[1mm]
-1&0&1&0&1\\[1mm]
\gamma^2-\gamma^{5/2}&\gamma-\sqrt{\gamma}&\gamma-1/\sqrt{\gamma}
&1&-\gamma^{3/2}\\[1mm]
\gamma^{-2}+i\gamma^{-5/2}&-\gamma^{-1}+i\gamma^{-1/2}
&-\gamma^{-1}-i\sqrt{\gamma}&1&-i\gamma^{-3/2}\\[1mm]
\gamma^{-2}-i\gamma^{-5/2}&-\gamma^{-1}-i\gamma^{-1/2}
&-\gamma^{-1}+i\sqrt{\gamma}&1&+i\gamma^{-3/2}
\end{array}\right).
\end{equation}
That is, $w_i\,M=\lambda_i\,w_i$ for $i=1,\ldots,5$.
Let
\begin{equation}\label{E4.07}
\begin{aligned}
C:=&(0,\,0,\,0,\,0,\,1)\cdot W^{-1}\\
=&\frac{1}{4\sqrt{5}}\left(\sqrt{\gamma}-1,\;2\sqrt{5},\;-1-\sqrt{\gamma},
\;1-i/\sqrt{\gamma},\;1+i/\sqrt{\gamma}\right).
\end{aligned}
\end{equation}
Summarizing the discussion, we get
$$\begin{aligned}
a_{n;1}&
=\sum_{j=1}^5C_j\,w_{j,1}\,\lambda_j^{n-1}\\
&=\frac{\gamma^2-\gamma^{3/2}}{4\sqrt{5}}\big(\gamma+\sqrt{\gamma}\big)^n
-\frac12
+\frac{\gamma^2+\gamma^{3/2}}{4\sqrt{5}}\big(\gamma-\sqrt{\gamma}\big)^n\\
&\quad
+\frac{-\gamma^{-2}-i\gamma^{-3/2}}{4\sqrt{5}}\big(-\gamma^{-1}+i/\sqrt{\gamma}\big)^n
+\frac{-\gamma^{-2}+i\gamma^{-3/2}}{4\sqrt{5}}\big(-\gamma^{-1}-i/\sqrt{\gamma}\big)^n.
\end{aligned}
$$
Hence, the proof of Proposition~\ref{P1.5} is completed.\eop

\begin{uremark}\rm
\label{R4.1}
By a standard procedure, the vector-valued linear recursion of order 1 in \equ{E4.02} can be transformed into a scalar-valued linear recursion of order 5:
\begin{equation}
\label{E4.08}
a_{n+1;1}=\sum_{j=0}^4\beta_j\,a_{n-j;1}
\end{equation}
for certain numbers $\beta_0,\ldots,\beta_4$. Using Proposition~\ref{P1.5}, we can compute
$$
\begin{array}{c||c|c|c|c|c|c|c|c|c|c|}
n&1&2&3&4&5&6&7&8&9&10\\\hline
a_{n;1}&0&0&1&4&12&36&105&304&880&2544\\
\end{array}
$$
Writing \equ{E4.08} for $n=5,6,7,8,9$ and inverting the matrix, we get
$$
\left(\begin{array}{c}
\beta_0\\\beta_1\\\beta_2\\\beta_3\\\beta_4\\
\end{array}\right)
=\left(\begin{array}{rrrrr}
12&4&1&0&0\\
36&12&4&1&0\\
105&36&12&4&1\\
304&105&36&12&4\\
880&304&105&36&12
\end{array}
\right)^{-1}
\left(\begin{array}{r}
36\\105\\304\\880\\2544\end{array}\right)
=
\left(\begin{array}{r}
3\\0\\0\\-3\\1\end{array}\right).$$

Hence $a_{n;1}$ is the solution of the recursion equation
\begin{equation}
\label{E4.09}
a_{n+1;1}=3\,a_{n;1}-3\,a_{n-3;1}+a_{n-4;1},\qquad n\geq5,\end{equation}
with initial values
\[a_{1;1}=a_{2;1}=0,\quad a_{3;1}=1,\quad a_{4;1}=4,\quad a_{5;1}=12.\tag*{$\Diamond$}\]
\end{uremark}
\subsection{The counting approach to probabilities: Proof of Theorem~\ref{T3}(i)}
\label{S4.2}
Recall $\alpha$, $\eta$ and $\hatf$ from Theorem~\ref{T3} and recall that $\P$ denotes the uniform spanning tree measure on $G^{H_3}$.
The aim of this section is to proceed similarly as in Section~\ref{S2.2} and Section~\ref{S3} to infer that the matrix $A$ of the determinantal point process $X$ is indeed given by $A_{k,l}=\hatf(l-k)$.

The statement will follow from the following three lemmas (recall that $\gamma=(\sqrt 5+1)/2$).
\begin{lemma}
\label{L4.2}
$\displaystyle \P[z_0\in T]=\hatf(0)=\frac{\gamma^{3/2}}{2\sqrt{5}}\approx0.460221$.
\end{lemma}
\begin{lemma}
\label{L4.3}
For $m\in\N$, we have
$\P[z_0,z_m\in T]=\hatf(0)^2-\hatf(m)^2$.
\end{lemma}
\begin{lemma}
\label{L4.4}For all $m\in\N$ and for $k=1,2$, we have
\begin{equation}
\label{E4.10}
P[z_0,z_k,z_{m+1}\in T]=
\det\left(
\begin{array}{ccc}
\hatf(0)           &\hatf(k)&\hatf(m+1)\\[1.5mm]
\hatf(k)    &       \hatf(0)&\hatf(m+1-k)    \\[1.5mm]
\hatf(m+1)&\hatf(m+1-k)&       \hatf(0)
\end{array}\right).
\end{equation}
\end{lemma}
\medskip

\textbf{Proof of Theorem~\ref{T3}(i).}\par
From Lemma~\ref{L4.2} we get $A_{0,0}=\hatf(0)>0$. From Lemma~\ref{L4.3} we get that
$|A_{0,m}|=\sqrt{\P[z_0,z_m\in T]-\hatf(0)^2}=|\hatf(m)|$ for $m\in \N$. As argued at the end of Section~\ref{S1.2}, we are free to choose the sign of $A_{0,1}$ and we make the choice $A_{0,1}=\hatf(1)$. Now we proceed by induction. Assume that we have shown already that $A_{0,n}=\hatf(n)$ for $n=0,\ldots,m$. Then, by Lemma~\ref{L4.4}, for $k=1,2$, the determinant on the right hand side of \equ{E4.10} equals
$$\det\left(
\begin{array}{ccc}
\hatf(0)           &\hatf(k)&A_{0,m+1}\\[1.5mm]
\hatf(k)    &       \hatf(0)&\hatf(m+1-k)    \\[1.5mm]
A_{0,m+1}&\hatf(m+1-k)&       \hatf(0)
\end{array}\right).
$$
Explicitly computing the determinant, we see that the sign of $A_{0,m+1}$ is determined by this equation and equals the sign of $\hatf(m+1)$ unless one of the other matrix entries is zero. Clearly, $\hatf(0)$, $\hatf(1)$ and $\hatf(2)$ are not zero. Now $\hatf(m)=0$ if and only if $\Re(\eta\,\alpha^m)=0$. However, in this case, since $\Im(\alpha)\neq0$, we have $\Re(\eta\,\alpha^{m-1})\neq0$. Concluding, we get $A_{0,m+1}=\hatf(m+1)$ from Lemma~\ref{L4.4} either using $k=1$ or $k=2$.

This finishes the proof of Theorem~\ref{T3}(i) subject to the Lemmas~\ref{L4.2}, \ref{L4.3} and \ref{L4.4}.\eop

It remains to prove Lemmas~\ref{L4.2}, \ref{L4.3} and \ref{L4.4}.\medskip

\textbf{Proof of Lemma~\ref{L4.2}.}\par

For $i,j=1,\ldots,5$ and $t^-\in\Gamma^-_{0;i}$, $t^+\in\Gamma^+_{0;j}$, define the set of bridges
\begin{equation}
\label{E4.11}
B_{i,j}:=\big\{B\subset \{h_1,h_2\}:\,t^-\cup B\cup t^+\in\SPT\big(G^{H_3}\big)\big\}.
\end{equation}
Let $N_{i,j}=\#B_{i,j}$ and compute that
\begin{equation}
\label{E4.12}
N=\left(\begin{array}{ccccc}
1&2&1&1&1\\
2&0&1&1&0\\
1&1&1&0&0\\
1&1&0&1&0\\
1&0&0&0&0
\end{array}\right)
.
\end{equation}

Recall $M$ and $M'$ from \equ{E4.01}
Arguing as in Section~\ref{S2.2}, we get
$$
\P[z_0\in T]=\frac{w_1M'N(w_1)^T}{\lambda_1\,w_1N(w_1)^T}=\frac{\gamma^{3/2}}{2\sqrt{5}}=\hatf(0).
$$
This finishes the proof of Lemma~\ref{L4.2}.\eop

\textbf{Proof of Lemma~\ref{L4.3}.}
Let $m\in\N$.
Arguing as in Section~\ref{S2.2}, we get
\begin{equation}
\label{E4.13}
\P[z_0,z_m\in T]
=\frac{w_1\,M'\,M^{m-1}\,N\,(w_1\,M')^T}{w_1\,M^{m+1}\,N\,w_1^T}.
\end{equation}
For a simpler representation, write
$w_1M'=\sum_{j=1}^5\mu_j\,w_j$
with
$(\mu_1,\ldots,\mu_5)=w_1M'W^{-1}$.
Define $\tilde\lambda_j:=\lambda_j/\lambda_1$ and
$$\varrho_j:=
\frac{\mu_j\,w_j\,N\,(w_1M')^T}
{\lambda_1\lambda_j\,w_1\,N\,(w_1)^T}.
$$
Then we have
$$
\P[z_0,z_m\in T]=\sum_{j=1}^5\varrho_j\,\tilde\lambda_j^m.
$$
An explicit computation yields
$$\begin{gathered}
\tilde\lambda_1=1,\qquad
\tilde\lambda_2=\gamma-\sqrt{\gamma},\qquad
\tilde\lambda_3=\gamma^3-2\gamma^{3/2},\\[2mm]
\tilde\lambda_4=(\gamma^{-1/2}-1)+i(\gamma^{1/2}-1),\qquad
\tilde\lambda_5=(\gamma^{-1/2}-1)-i(\gamma^{1/2}-1),
\end{gathered}
$$
and
$$
\varrho_1=\frac{\gamma^3}{20},
\quad \varrho_2=-\frac{1}{4\sqrt{5}},
\quad \varrho_3=0,\quad \varrho_4=\frac{i-2}{40},
\quad\varrho_5=\frac{-i-2}{40}.
$$
Note that $\tilde\lambda_3=\tilde\lambda_2^2$ and that $|\tilde\lambda_4|=|\tilde\lambda_5|=\lambda_2$. Furthermore, $\varrho_2^2=4\varrho_4\varrho_5$.
Note that $\eta=\sqrt{-\varrho_5}$ and that $\alpha=\sqrt{\tilde\lambda_5}$.
A simple computation yields $\hatf(0)^2-\hatf(m)^2=\sum_{j=1}^5\varrho_j\,\tilde\lambda_j^m$.
This completes the proof of Lemma~\ref{L4.3}.\eop

\textbf{Proof of Lemma~\ref{L4.4}.}
Iterating the argument in Section~\ref{S2.2}, we get the probability that three edges $0<k<m+1$ are in the spanning tree
\begin{equation}
\label{E4.14}
\P\big[z_0,\,z_k,\,z_{m+1}\in T\big]=\sum_{j=1}^5\varrho^{(k)}_j\tilde\lambda_j^m
\end{equation}
with $$\varrho^{(k)}_j:=
\frac{(w_1M'M^{k-1}M'W^{-1})_j\;\mu_j\,w_j\,N\,(w_1M')^T}
{\lambda_1^2\lambda_j\,w_1\,N\,(w_1)^T}.
$$
Explicit computations for $k=1$ and $k=2$ yield
$$
\varrho^{(1)}_1=\frac{\gamma^{5/2}}{40},\quad
\varrho^{(1)}_2=\frac{1-\gamma^{3/2}}{4\sqrt5},\quad
\varrho^{(1)}_3=0,\quad
\varrho^{(1)}_4=\frac{\sqrt{\gamma}-\sqrt5+i/\sqrt{\gamma}}{80},\quad
\varrho^{(1)}_5=\overline{\varrho^{(1)}_4}.
$$
and
$$\begin{gathered}
\varrho^{(2)}_1=-\frac{\gamma^{5/2}}{20}+\frac{\gamma^2}{10},\quad
\varrho^{(2)}_2=-\frac{\gamma^{-1}+\gamma^{-1/2}}{4\sqrt5},\quad
\varrho^{(2)}_3=0,\\
\varrho^{(2)}_4=\frac{\gamma^{-3/2}+i\gamma^{3/2}+2-i}{40},\quad
\varrho^{(2)}_5=\overline{\varrho^{(2)}_4}.
\end{gathered}
$$Recalling that $\alpha^2=\tilde\lambda_5$, $\bar\alpha^2=\tilde\lambda_4$ and $\alpha\bar \alpha=\tilde\lambda_2$, a simple computation using \equ{E4.14} yields the claim of Lemma~\ref{L4.4}.\eop

\subsection{Proof of Theorem~\ref{T3}(ii)}
\label{S4.3}

Let $f$ be defined by \equ{E1.22} and for $m\in\N_0$, let
$$\tilde f(m):=
\int_0^1e^{2\pi i mx}f(x)\,dx=
\int_0^1\frac{e^{2\pi i mx}}{4+4\cos(2\pi x)+2\cos(4\pi x)}\,dx.
$$
By the Fourier inversion formula, it is enough to check that $\tilde f(m)=\hatf(m)$ for $m\in\N_0$.

Using the substitution $y=e^{2\pi i x}$, we get

$$\begin{aligned}
&\tilde f(m)=\int_0^1\frac{e^{2\pi i mx}}
{4+2e^{2\pi i x}+2e^{-2\pi i x}+e^{4\pi i x}+e^{-4\pi i x}}\,dx\\
&\hspace*{12em}=\frac{1}{2\pi i}\oint\frac{y^m}{4+2y+2y^{-1}+y^2+y^{-2}}\,\frac{dy}{y}\\
&\hspace*{12em}=\frac{1}{2\pi i}\oint\frac{y^{m+1}}{y^4+2y^3+4y^2+2y+1}\,dy,
\end{aligned}
$$
where $\oint$ denotes the (anti-clockwise) curve integral along the unit sphere in the complex plane. The polynomial in the denominator can be decomposed into linear factors
$$g(y):=y^4+2y^3+4y^2+2y+1=(y-y_1)(y-y_2)(y-y_3)(y-y_4),$$
with
$$
y_{1,2}=\frac{\gamma^{-3/2}-1}{2}\pm i\frac{\gamma^{3/2}-1}{2}\mbu
y_{3,4}=-\frac{1+\gamma^{-3/2}}{2}\pm i\frac{1+\gamma^{3/2}}{2}.
$$
Note that $|y_{1,2}|=\sqrt{\gamma-\sqrt{\gamma}}<1$ and
$|y_{3,4}|=\sqrt{\gamma+\sqrt{\gamma}}>1$. Hence $y_1$ and $y_2$ are in the domain of integration and residue calculus yields
$$\begin{aligned}
\tilde f(m)
&=\frac{1}{2\pi i}\oint\frac{y^{m+1}}{g(y)}\,dy
=\frac{y_1^{m+1}}{g'(y_1)}+\frac{y_2^{m+1}}{g'(y_2)}\\
&=\frac{y_1^{m+1}}{4y_1^3+6y_1^2+8y_1+2}+\frac{y_2^{m+1}}{4y_2^3+6y_2^2+8y_2+2}\\
&=\eta\,\alpha^m+\bar\eta\,\bar\alpha^m\;=\;\hatf(m).
\end{aligned}
$$
This finishes the proof of Theorem~\ref{T3}(ii).\eop
\Section{The enhanced Helix-3-graph}
\label{S5}
The proof of Theorem~\ref{T4} is similar to that of Theorem~\ref{T3} but the actual computations are annoyingly tedious, although straightforward. Hence, we only give the main steps of the proof.

The idea is to employ the same counting scheme as in Section~\ref{S4}. Define the sets $F_{i,j}$ and $F'_{i,j}$ and matrices $M$ and $M'$ as in Section~\ref{S4}.
\begin{figure}
\label{F5.1}
\begin{picture}(400,140)(-50,0)
\put(20,02){\includegraphics[scale=1.2]{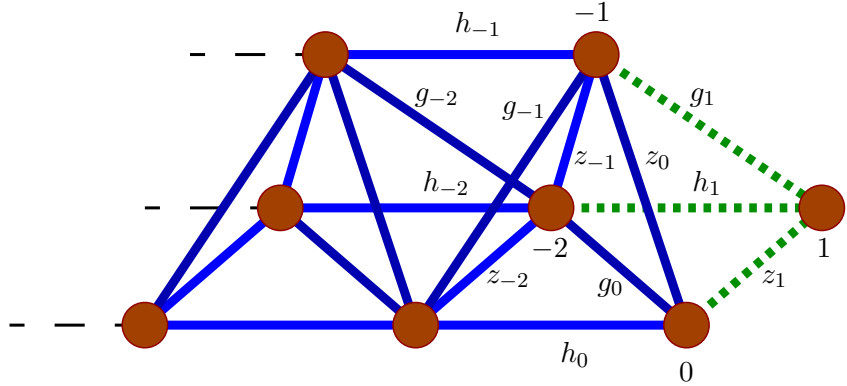}}
\put(275,06){$0$}
\put(219,53){$-2$}
\put(327,53){$1$}
\put(235,142){$-1$}
\put(230,13){$h_{0}$}
\put(175,112){$g_{-2}$}
\put(208,107){$g_{-1}$}
\put(244,40){$g_{0}$}
\put(279,112){$g_{1}$}
\put(178,78){$h_{-2}$}
\put(280,78){$h_{1}$}
\put(190,138){$h_{-1}$}
\put(202,43){$z_{-2}$}
\put(306,43){$z_{1}$}
\put(262,88){$z_{0}$}
\put(235,88){$z_{-1}$}
\end{picture}
\caption[b]{Adding some of the edges $\{g_1,h_1,z_1\}$ changes how $0,-1,-2$ are connected.}
\end{figure}
\begin{table}[ht]
$$
\def\arraycolsep{1.5mm}
\begin{array}{c|c|c|c|c|c|}
i\;\backslash\; j&1&2&3&4&5\\\hline
1&\{\{z_1\},\{h_1\},\{g_1\}\}&\{\emptyset\}&\emptyset&\emptyset&\emptyset\\[1.5mm]
2&\{\{g_1,z_1\},\{h_1,z_1\}\}&\emptyset&\{\{z_1\}\}&\{\{g_1\},\{h_1\}\}&\{\emptyset\}\\[1.5mm]
3&\{\{h_1,z_1\},\{g_1,h_1\}\}&\{\{h_1\}\}&\emptyset&\emptyset&\emptyset\\[1.5mm]
4&\{\{g_1,h_1\},\{g_1,z_1\}\}&\emptyset&\{\{h_1\},\{z_1\}\}&\{\{g_1\}\}&\{\emptyset\}\\[1.5mm]
5&\{\{g_1,h_1,z_1\}\}&\emptyset&\{\{h_1,z_1\}\}&\{\{g_1,h_1\}\}&\{\{h_1\}\}\\\hline
\end{array}
$$
\caption{The sets $F_{i,j}$ for the enhanced helix-3-graph}
\label{Tab5.1}
\end{table}
Explicit counting (see Table~\ref{Tab5.1} for $F_{i,j}$) yields
\begin{equation}
\label{E5.01}
M=\left(
\begin{array}{ccccc}
c+d+1& 1& 0& 0& 0\\
c(d+1)& 0& c& d+1& 1\\
c+d& 1& 0& 0& 0\\
(c+1)d&0& c+1& d& 1\\
cd& 0& c& d& 1
\end{array}
\right)\!\!\!\!\!\!\mbu\!\!\!\!
M'=\left(
\begin{array}{ccccc}
c& 0& 0& 0& 0\\
c(d+1)& 0& c& 0& 0\\
c& 0& 0& 0& 0\\
cd& 0& c& 0& 0\\
cd& 0& c& 0& 0
\end{array}
\right)
.
\end{equation}
The characteristic polynomial of $M$ is
\begin{equation}
\label{E5.02}
\begin{aligned}
\chi_M(x)&=x^5-(2d+c+2)x^4+(1-c+2d+d^2)x^3+(c-1-2d-d^2)x^2\\&\ +(c+2+2d)x-1\\
&=(x-1)\big(x^4-(2d+c+1)x^3-(2c-d^2)x^2-(2d+c+1)x+1\big).
\end{aligned}
\end{equation}
Let $$g(x)=c^{-1}\big[x^4+(d+2)x^3+(c+2d+3)x^2+(d+2)x+1]$$
and let $x_k$, $k=1,2,3,4$ be the roots of $g$. That is
$$x_{1}=\frac{1}{4}\big(-d-2+\sqrt {{d}^{2}-4\,d-4\,c}\big)+\frac14\sqrt { \left( -d
-2+\sqrt {{d}^{2}-4d-4\,c} \right) ^{2}-16},
$$
$$x_{2}=\frac{1}{4}\big(-d-2-\sqrt {{d}^{2}-4\,d-4\,c}\big)+\frac14\sqrt { \left( -d
-2-\sqrt {{d}^{2}-4d-4\,c} \right) ^{2}-16},
$$
$$x_{3}=\frac{1}{4}\big(-d-2+\sqrt {{d}^{2}-4\,d-4\,c}\big)-\frac14\sqrt { \left( -d
-2+\sqrt {{d}^{2}-4d-4\,c} \right) ^{2}-16},
$$
and
$$x_{4}=\frac{1}{4}\big(-d-2-\sqrt {{d}^{2}-4\,d-4\,c}\big)-\frac14\sqrt { \left( -d
-2-\sqrt {{d}^{2}-4d-4\,c} \right) ^{2}-16}
$$
Is is straightforward to verify that the roots of $\chi_M$ are
\begin{equation}
\label{E5.03}
\begin{aligned}
\lambda_1&=1/(x_1x_2),\\
\lambda_2&=1,\\
\lambda_3&=x_1x_2,\\
\lambda_4&=x_1/x_2,\\
\lambda_5&=x_2/x_1.
\end{aligned}
\end{equation}
Let $N_{i,j}=\weight(B_{i,j})$ as in Section~\ref{S4} and compute that
\begin{equation}
\label{E5.04}
N=\left(\begin{array}{ccccc}
1&d+2&1&d+1&d+1\\
d+2&0&d+1&d+1&0\\
1&d+1&1&d&d\\
d+1&d+1&d&2d+1&d\\
d+1&0&d&d&0
\end{array}\right)
.
\end{equation}
Define $\varrho_j$, $j=1,\ldots,5$ as in the proof of Lemma~\ref{L4.3} and note that $\varrho_3=0$ and $\varrho_2^2=4\varrho_4\varrho_5$. An explicit calculation yields
$$\varrho_4=-\left(\frac{x_1}{g'(x_1)}\right)^2\mbu\varrho_5=-\left(\frac{x_2}{g'(x_2)}\right)^2.$$
 Arguing as in the proof of Lemma~\ref{L4.3} we get that $\hatf$ is given for $m\in\N_0$ by
$$\hatf(m)=\sqrt{\varrho_2(\lambda_2/\lambda_1)^m+\varrho_4(\lambda_4/\lambda_1)^m+\varrho_5(\lambda_5/\lambda_1)^m}=\sqrt{-\varrho_4}\;x_1^m+\sqrt{-\varrho_5}\;x_2^m.$$
Now proceed as in the proof of Lemma~\ref{L4.4} to compute that $\hatf$ is in fact the inverse Fourier transform of $f$ from \equ{E1.25}:
$$\int_0^1e^{2\pi imx}f(x)\,dx=\frac{1}{2\pi i}\oint\frac{y^{m+1}}{g(y)}\,dy
=\frac{x_1^{m+1}}{g'(x_1)}+\frac{x_2^{m+1}}{g'(x_2)}=\hatf(m).
$$
\Section{Helix-3-graph: The Markov chain approach.}
\label{S6}

Following H\"{a}ggstr\"{o}m \cite{Haggstrom1994}, we consider a subset $E'\subset E^{H_3}$ of the edges of $G^{H_3}$ as an element of the product space $\Pi:=(\{0,1\}\times\{0,1\})^\Z$. Denote by $(h,z)=(\bar h_n,\bar z_n)_{n\in\Z}$ a generic element of that space that corresponds to $E'$ via $h_n\in E'$ iff $\bar h_n=1$ and $z_n\in E'$ iff $\bar z_n=1$. Denote by $\vartheta:(\{0,1\}\times\{0,1\})^\Z\to(\{0,1\}\times\{0,1\})^\Z$, $(\bar h_n,\bar z_n)_{n\in \Z}\mapsto (\bar h_{n+1},\bar z_{n+1})_{n\in\Z}$ the shift operator. Then the uniform spanning tree measure $\P$ is invariant for $\vartheta$. That is, $(\P,\vartheta)$ is a measure preserving dynamical system. $\P$ is concentrated on the set $\SPT\big(G^{H_3}\big)$ of configurations that are spanning trees. $\SPT\big(G^{H_3}\big)$ can be characterized as subset of $\Pi$ where a certain translation invariant dictionary of finite letter words (with alphabet $\{0,1\}\times\{0,1\}$) are forbidden. For such a situation, as pointed out in \cite[Theorems 2.4 and 2.5]{Haggstrom1994} (for proofs, see \cite{BurtonSteif1994,Haggstrom1995,Parry1964}) the uniform distribution on the allowed configurations $\SPT\big(G^{H_3}\big)$ can be characterized as the measure $\P$ concentrated on $\Pi$ such that $(\P,\vartheta)$ is a measure preserving dynamical system with \emph{maximal entropy}. Furthermore, $(\P,\vartheta)$ is a stationary Markov chain if all forbidden words have length at most 2. Since words of length 2 do not suffice to decide whether a given configuration is a spanning tree or not (in fact, arbitrarily long words are needed), we will develop a more subtle encoding of the spanning tree that yields the Markov property.

Let $t\in\SPT\big(G^{H_3}\big)$ be a spanning tree on $G^{H_3}$. For each $n$, we have $t\rest{\{n,n-1,n-2,\ldots\}}\in\Gamma^-_{n;i}$ for exactly one $i=i(t,n)=1,\ldots,5$. Define $\vp_n(t)=(\1_t(h_n),\1_t(z_n),i(t,n))$. Note that $\vp_n(t)$ takes values in the set
$$\Delta:=\{\delta_k:\,k=1,\ldots,11\},$$
where the symbols $(\delta_k)$ are defined by the following table.
\begin{equation}
\label{E6.01}
\begin{array}{c|c}
k&\delta_k\\\hline
1&(0,1,1)\\
2&(1,0,1)\\
3&(1,1,1)\\[2mm]
4&(0,0,2)\\
5&(1,0,2)\\
&
\end{array}
\quad
\begin{array}{r|c}
k&\delta_k\\\hline
6&(0,1,3)\\
7&(1,0,3)\\
8&(1,1,3)\\[2mm]
9&(1,0,4)\\[2mm]
10&(0,0,5)\\
11&(1,0,5)
\end{array}
\end{equation}

We group the elements of $\Delta$ according to their third entry:
$$
\Delta_1:=\{\delta_1,\delta_2,\delta_3\},\quad
\Delta_2:=\{\delta_4,\delta_5\},\quad
\Delta_3:=\{\delta_6,\delta_7,\delta_8\},\quad
\Delta_4:=\{\delta_9\},\quad
\Delta_5:=\{\delta_{10},\delta_{11}\}.
$$
A symbol $\delta\in\Delta_1$ cannot be followed by $\delta_3$ or $\delta_8$ as this would create a cycle. If followed by $(h,z)=(0,0)$, the new symbol is necessarily in $\Delta_2$, hence $\delta_{10}$ is forbidden. Going through all possibilities, we get the following list of possible successors of $\delta\in\Delta_k$.
\begin{equation}
\label{E6.02}
\begin{array}{c|l}
k&\mbox{allowed successors}\\\hline
1&\delta_1,\,\delta_2,\,\delta_4\\
2&\delta_3,\,\delta_6,\,\delta_9,\,\delta_{10}\\
3&\delta_3,\,\delta_5\\
4&\delta_6,\,\delta_7,\,\delta_{10}\\
5&\delta_8,\,\delta_{11}
\end{array}
\end{equation}

Let $D\subset\Delta^\Z$ denote the subset where all two-letter words not in the above list are forbidden. Then $\vp=(\vp_n)_{n\in\Z}:\SPT\big(G^{H_3}\big)\to D$ is a bijection.
Furthermore, the image measure $\P':=\P\circ \vp^{-1}$ on $D$ is the unique maximizer of the entropy among all stationary measures on $D$ and the canonical process $Y=(Y_n)_{n\in\Z}$ is a stationary Markov chain under $\P'$. Since any stationary measure on $D$ defines a Markov chain with state space $\Delta$, we can characterize $\P'$ as the distribution that maximizes the entropy among all Markov chains on $\Delta$ with the allowed transitions given in the table \equ{E6.02}. In order to compute this entropy, let $R=(R_{i,j})_{i,j=1,\ldots,11}$ denote the transition matrix of the Markov chain (under $\P'$) and denote by $\pi^R$ the invariant distribution (on $\Delta$). Recall that the entropy is
\begin{equation}
\label{E6.03}
H^R=-\sum_{i,j=1}^{11}\pi^R_iR_{i,j}\log(R_{i,j}).
\end{equation}
Taking into account the natural symmetries and recalling that the uniform distribution on a finite set maximizes the entropy, we get that an $R$ that maximizes $H^R$ has to be of the form
\begin{equation}
\label{E6.04}
\newcommand{\vsc}[1]{\scriptstyle{#1}}
\def\arraycolsep{2pt}
R=\left(\begin{array}{ccccccccccc}
\scriptstyle{\bar R_{1,1}/2}&\vsc{\bar R_{1,1}/2}&\vsc{0}&\vsc{1-\bar R_{1,1}}&\vsc{0}&\vsc{0}&\vsc{0}&\vsc{0}&\vsc{0}&\vsc{0}&\vsc{0}\\
\vsc{\bar R_{1,1}/2}&\vsc{\bar R_{1,1}/2}&\vsc{0}&\vsc{1-\bar R_{1,1}}&\vsc{0}&\vsc{0}&\vsc{0}&\vsc{0}&\vsc{0}&\vsc{0}&\vsc{0}\\
\vsc{\bar R_{1,1}/2}&\vsc{\bar R_{1,1}/2}&\vsc{0}&\vsc{1-\bar R_{1,1}}&\vsc{0}&\vsc{0}&\vsc{0}&\vsc{0}&\vsc{0}&\vsc{0}&\vsc{0}\\
\vsc{0}&\vsc{0}&\vsc{\bar R_{2,1}}&\vsc{0}&\vsc{0}&\vsc{\bar R_{2,3}}&\vsc{0}&\vsc{0}&\vsc{\bar R_{2,4}}&\vsc{1-\bar R_{2,1}-\bar R_{2,3}-\bar R_{2,4}}&\vsc{0}\\
\vsc{0}&\vsc{0}&\vsc{\bar R_{2,1}}&\vsc{0}&\vsc{0}&\vsc{\bar R_{2,3}}&\vsc{0}&\vsc{0}&\vsc{\bar R_{2,4}}&\vsc{1-\bar R_{2,1}-\bar R_{2,3}-\bar R_{2,4}}&\vsc{0}\\
\vsc{0}&\vsc{0}&\vsc{\bar R_{3,1}}&\vsc{0}&\vsc{1-\bar R_{3,1}}&\vsc{0}&\vsc{0}&\vsc{0}&\vsc{0}&\vsc{0}&\vsc{0}\\
\vsc{0}&\vsc{0}&\vsc{\bar R_{3,1}}&\vsc{0}&\vsc{1-\bar R_{3,1}}&\vsc{0}&\vsc{0}&\vsc{0}&\vsc{0}&\vsc{0}&\vsc{0}\\
\vsc{0}&\vsc{0}&\vsc{\bar R_{3,1}}&\vsc{0}&\vsc{1-\bar R_{3,1}}&\vsc{0}&\vsc{0}&\vsc{0}&\vsc{0}&\vsc{0}&\vsc{0}\\
\vsc{0}&\vsc{0}&\vsc{0}&\vsc{0}&\vsc{0}&\vsc{\bar R_{4,3}/2}&\vsc{\bar R_{4,3}/2}&\vsc{0}&\vsc{0}&\vsc{1-\bar R_{4,3}}&\vsc{0}\\
\vsc{0}&\vsc{0}&\vsc{0}&\vsc{0}&\vsc{0}&\vsc{0}&\vsc{0}&\vsc{\bar R_{5,3}}&\vsc{0}&\vsc{0}&\vsc{1-\bar R_{5,3}}\\
\vsc{0}&\vsc{0}&\vsc{0}&\vsc{0}&\vsc{0}&\vsc{0}&\vsc{0}&\vsc{\bar R_{5,3}}&\vsc{0}&\vsc{0}&\vsc{1-\bar R_{5,3}}
\end{array}
\right)
\end{equation}
where $\bar R_{1,1},\bar R_{2,1},\bar R_{2,3},\bar R_{2,4},\bar R_{3,1},\bar R_{4,3},\bar R_{5,3}$ are the seven free parameters in the problem. It seems hopeless to solve the entropy maximizing problem for these seven parameters analytically. Numerically, however, the problem is rather easy to compute.
\ignore{
to get a closed form of the parameters that solve the maximizing problem for \equ{E6.03} but with a numerical scheme using MAPLE we could compute the values (in fact, with 100 digits precision):
\begin{equation}
\label{E6.05}
\begin{aligned}
\bar R_{1,1}&=0.69202867847165176793\\
\bar R_{2,1}&=0.38875672696616535394\\
\bar R_{2,3}&=0.25424132496148527486\\
\bar R_{2,4}&=0.22248654606766929211\\
\bar R_{3,1}&=0.52908551363574612516\\
\bar R_{4,3}&=0.79079967411811755990\\
\bar R_{5,3}&=0.65398566076417411609
\end{aligned}
\end{equation}
and
\begin{equation}
\label{E6.06}
\;\;H^R=1.06127506190503565203.
\end{equation}
}
The precision is easily made good enough for all applications. For example, we can draw random samples of uniform spanning trees on $G^{H_3}_{0,n}$ for arbitrarily large $n$ by simulating the Markov chain. See Figure~\ref{F5.1}.

\begin{figure}[h]
\centerline{\includegraphics[width=0.9\textwidth]{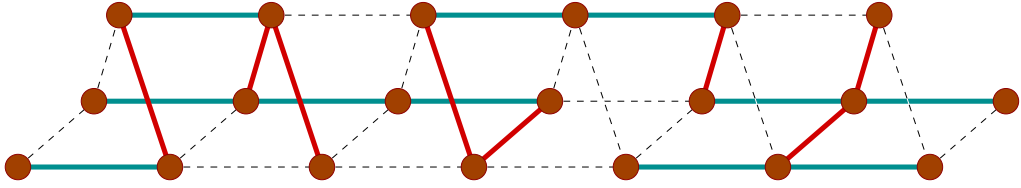}}
\caption[b]{A random spanning tree on 20 points. Circular edges $(z_n)$ are colored red and horizontal edges $(h_n)$ are colored blue.}
\label{F6.1}
\end{figure}

Now we present a way to obtain the exact solution of the maximizing problem for \equ{E6.03} in a more subtle way. Consider the projection $\mathrm{pr}:\Delta\to\{1,\ldots,5\}$, $(h,z,i)\mapsto i$ and note that $\mathrm{pr}^{-1}(\{i\})=\Delta_i$. Note that due to the symmetries, also $\bar Y_n:=\mathrm{pr}(Y_n)$ is a Markov chain on $\{1,\ldots,5\}$ with transition matrix
\begin{equation}
\label{E6.07}
\bar R=\left(\begin{array}{ccccc}
\bar R_{1,1}&1-\bar R_{1,1}&0&0&0\\
\bar R_{2,1}&0&\bar R_{2,3}&\bar R_{2,4}&1-\bar R_{2,1}-\bar R_{2,3}-\bar R_{2,4}\\
\bar R_{3,1}&1-\bar R_{3,1}&0&0&0\\
0&0&\bar R_{4,3}&0&1-\bar R_{4,3}\\
0&0&\bar R_{5,3}&0&1-\bar R_{5,3}
\end{array}
\right).
\end{equation}
Denote by $\pi^{\bar R}$ the invariant distribution of this chain and write
\begin{equation}
\label{E6.08}
\tilde \pi^{\bar R}_{i,j}=\pi^{\bar R}_i\bar R_{i,j}\mf i,j=1,\ldots,5,
\end{equation}
for the invariant distribution of the bivariate chain $(\bar Y_n,\bar Y_{n+1})_{n\in\Z}$.
Note that $\pi^{\bar R}_i$ and $\tilde\pi^{\bar R}_{i,j}$ are quantities that we can read off directly from the uniform spanning tree measure:

\begin{equation}
\label{E6.09}
\begin{aligned}
\pi^{\bar R}_i&=P\Big(T\rest{\{-\infty,\ldots,0\}}\in\Gamma^-_{0;i}\Big),\\
\tilde\pi^{\bar R}_{i,j}&=P\Big(T\rest{\{-\infty,\ldots,0\}}\in\Gamma^-_{0;i},
\;T\rest{\{-\infty,\ldots,1\}}\in\Gamma^-_{1;j}\Big).
\end{aligned}
\end{equation}
The counting scheme from Section~\ref{S4.2} yields
\begin{equation}
\label{E6.10}
\pi^{\bar R}_i=\frac{w_{1,i}\,(N\,w_{1}^T)_i}{w_1\,N\,w_1^T}.
\end{equation}
A direct computation yields
\begin{equation}
\label{E6.11}
\begin{aligned}
\pi^{\bar R}&=\frac{1}{4\sqrt{5}}\Big(\gamma^2+
\gamma^{3/2},\,
\gamma^{-1/2}+
2\gamma^{-1},\,
2\gamma^{-1},\,
2\gamma^{-1}-\gamma^{-1/2},\,
\gamma^2-\gamma^{3/2}\Big)\\
&=\big(0.522816,\, 0.226091,\, 0.1381966,\, 0.050302, \,0.06259458\big).
\end{aligned}
\end{equation}
Furthermore, we get
\begin{equation}
\label{E6.12}
\tilde\pi^{\bar R}_{i,j}=\frac{w_{1,i}\,M_{i,j}\,(N\,w_1^T)_j}
{w_1\,M\,N\,w_1^T}.
\end{equation}
Recall that $w_1\,M=\lambda_1\,w_1$ to simplify the expression and get
\begin{equation}
\label{E6.13}
\bar R_{i,j}=\frac{\tilde\pi^{\bar R}_{i,j}}{\pi^{\bar R}_i}=
\frac{w_{1,i}\,M_{i,j}\,(N\,w_{1}^T)_j}
{\lambda_1w_{1,i}\,(N\,w_{1}^T)_i}.
\end{equation}
Hence the non-zero entries of $\bar R$ are
\begin{equation}
\label{E6.14}
\begin{array}{cll}
\bar R_{1,1}&=2(\gamma-\sqrt{\gamma})       &=0.6920287,\\[2mm]
\bar R_{1,2}&=1-2(\gamma-\sqrt{\gamma})      &=0.3079713,\\[2mm]
\bar R_{2,1}&=\frac{\sqrt{\gamma}}{2+\sqrt{\gamma}}&=0.3887567,\\[2mm]
\bar R_{2,3}&=\frac{\gamma-1/\sqrt{\gamma}}{2+\sqrt{\gamma}}&=0.2542413,\\[2mm]
\bar R_{2,4}&=\frac{2-\sqrt{\gamma}}{2+\sqrt{\gamma}}&=0.2224865,\\[2mm]
\bar R_{2,5}&=\frac{\sqrt{\gamma}+1/\sqrt{\gamma}-\gamma}{2+\sqrt{\gamma}}&=0.1345154,\\[2mm]
\bar R_{3,1}&=\frac{\gamma^{3/2}-1}{2}&=0.5290855,\\[2mm]
\bar R_{3,2}&=\frac{3-\gamma^{3/2}}{2}&=0.4709145,\\[2mm]
\bar R_{4,3}&=\frac{2\sqrt{5}\,\gamma-2\gamma^{5/2}}{2-\sqrt{\gamma}}&=0.7907997,\\[2mm]
\bar R_{4,5}&=\frac{2-\sqrt{\gamma}-2\sqrt{5}\,\gamma+2\gamma^{5/2}}{2-\sqrt{\gamma}}&=0.2092003,\\[2mm]
\bar R_{5,3}&=1-\frac{1}{\gamma+\sqrt{\gamma}}&=0.6539857,\\[2mm]
\bar R_{5,5}&=\frac{1}{\gamma+\sqrt{\gamma}}&=0.3460143.
\end{array}
\end{equation}
The invariant distribution of $R$ (recall \equ{E6.04}) is
\begin{equation}
\label{E6.15}
\begin{aligned}
\pi^{R}&=\frac{1}{4\sqrt{5}}\Big(\gamma,\gamma,
\gamma^{3/2}-\gamma^{-1},
\gamma^{3/2}-\gamma^{-1},
3\gamma^{-1}-\gamma^{1/2},
2\gamma-2\gamma^{1/2},\\
&\qquad\qquad\sqrt{5}-\gamma^{3/2},
\sqrt{5}\,\gamma^{3/2}-\gamma^3,
2\gamma^{-1}-\gamma^{-1/2},
\sqrt{5}\,\gamma^{3/2}-\gamma^3,
\gamma^4-2\gamma^{5/2}\Big)\\
&=\big(0.180902, 0.180902, 0.161012, 0.161012, 0.0650785, 0.0773725,\\&
 \phantom{=\;(}  0.019890, 0.040929, 0.050302, 0.040929, 0.0216675\big).
\end{aligned}
\end{equation}

The entropy can be computed using \equ{E6.03} or can be deduced simply from the fact that the entropy of a uniform distribution is the logarithm of the size of the state space. Either way we get (recall \equ{E4.05})
\begin{equation}
\label{E6.16}
\begin{aligned}
H^R&=\log(\lambda_1)=\log\big(\gamma+\sqrt{\gamma}\,\big)
=1.061\ldots
\end{aligned}
\end{equation}
We have thus shown the following theorem.
\begin{theorem}
\label{T5}
The process $Y=(Y_n)_{n\in\Z}$ is a stationary Markov chain with transition matrix $R$ given by \equ{E6.04} and \equ{E6.14}. The invariant distribution is given by \equ{E6.15}, the entropy by \equ{E6.16}.
\end{theorem}
The Markov chain can be used for simulations as well as for explicit computations. For example, we compute (for $m\in\N$)
\begin{equation}
\label{E6.17}
\begin{aligned}
\P[z_0\in T]
&=\sum_{i\in\{1,3,6,8\}}\pi^R_i=\frac{\gamma^{3/2}}{2\sqrt{5}}=0.460221\\
\P[z_0,z_1\in T]
&=\sum_{i\in\{1,3,6,8\}}\pi^R_i\,(R_{i,1}+R_{i,3})=\frac{\gamma}{4\sqrt{5}}=0.180902\\
\P[z_0,z_1,\ldots,z_m\in T]
&=\sum_{i\in\{1,3,6,8\}}\pi^R_i\,(R_{i,1}+R_{i,3})R_{1,1}^{m-1}
=\frac{\gamma}{4\sqrt{5}}(\gamma-\sqrt{\gamma})^{m-1}\\
&=\frac{\gamma^2+\gamma^{3/2}}{4\sqrt{5}}\big(\gamma-\sqrt{\gamma}\big)^{m}
=0.52281560\cdot0.3460143^m.
\end{aligned}
\end{equation}
Note that, in fact, $R_{1,1}=R_{3,1}$ which yields the simple form in the last equation.

\bibliographystyle{plain}
\def\cprime{$'$}

\end{document}